\documentclass[reqno]{amsart}

\usepackage{a4wide}
\usepackage{hyperref}
\usepackage{cleveref}
\usepackage{newtxtext}
\usepackage{newtxmath}
\usepackage{tikz}
\usetikzlibrary{cd}
\usepackage[mathscr]{euscript}
\usepackage{mathtools}
\usepackage{subcaption}

\usepackage[T1]{fontenc}

\newcommand{\graf}{\mathsf{\Gamma}}

\newcommand{\sfC}{\mathsf{C}}

\newcommand{\sfT}{\mathsf{T}}

\newcommand{\sfa}{\mathsf{a}}
\newcommand{\sfb}{\mathsf{b}}

\newcommand{\sfe}{\mathsf{e}}

\newcommand{\sfh}{\mathsf{h}}

\newcommand{\sfv}{\mathsf{v}}
\newcommand{\sfw}{\mathsf{w}}

\newcommand{\bbF}{\mathbb{F}}

\newcommand{\bbR}{\mathbb{R}}
\newcommand{\bbS}{\mathbb{S}}
\newcommand{\bbZ}{\mathbb{Z}}

\newcommand{\cH}{\mathcal{H}}

\newcommand{\cT}{\mathcal{T}}
\newcommand{\cV}{\mathcal{V}}

\newcommand{\bfa}{\mathbf{a}}
\newcommand{\bfb}{\mathbf{b}}
\newcommand{\bfe}{\mathbf{e}}
\newcommand{\bfi}{\mathbf{i}}
\newcommand{\bfk}{\mathbf{k}}
\newcommand{\bfs}{\mathbf{s}}
\newcommand{\bft}{\mathbf{t}}

\newcommand{\bfx}{\mathbf{x}}

\newcommand{\ess}{\mathsf{ess}}

\newcommand{\conv}{\mathop{\scalebox{1.5}{\raisebox{-0.2ex}{\(\ast\)}}}}

\newcommand{\grape}[2][0]{
\draw[thick, fill] (#2.center) circle(2pt) -- +({30+#1}:0.5) +(0,0) -- +({-30+#1}:0.5);
\draw[thick] (#2.center) +({30+#1}:0.5) arc ({120+#1}:{-120+#1}:{0.5 / sqrt(3)});
}

\newcommand\loops{\ell}
\newcommand\TC{\sfC^{\mathsf{top}}}

\newcommand{\field}{\bbF}
\newcommand{\ring}{R}

\newcommand\DGA{\(\mathsf{DGAlg}\)}
\newcommand\DGCA{\(\mathsf{DGCoAlg}\)}

\DeclareMathOperator{\val}{val}
\DeclareMathOperator{\image}{Im}

\DeclareMathOperator{\pr}{Proj}
\DeclareMathOperator{\Star}{Star}
\DeclareMathOperator{\Loop}{Loop}
\DeclareMathOperator{\SL}{SL}

\numberwithin{equation}{section}

\theoremstyle{plain}
\newtheorem{theorem}{Theorem}[section]
\newtheorem{lemma}[theorem]{Lemma}
\newtheorem{corollary}[theorem]{Corollary}
\newtheorem{proposition}[theorem]{Proposition}

\theoremstyle{definition}
\newtheorem{definition}[theorem]{Definition}
\newtheorem{assumption}[theorem]{Assumption}

\newtheorem{example}[theorem]{Example}

\theoremstyle{remark}
\newtheorem{remark}[theorem]{Remark}

\AddToHook{env/lemma/begin}{\crefalias{theorem}{lemma}}
\AddToHook{env/corollary/begin}{\crefalias{theorem}{corollary}}
\AddToHook{env/proposition/begin}{\crefalias{theorem}{proposition}}
\AddToHook{env/claim/begin}{\crefalias{theorem}{claim}}
\AddToHook{env/definition/begin}{\crefalias{theorem}{definition}}
\AddToHook{env/assumption/begin}{\crefalias{theorem}{assumption}}
\AddToHook{env/convention/begin}{\crefalias{theorem}{convention}}
\AddToHook{env/observation/begin}{\crefalias{theorem}{observation}}
\AddToHook{env/example/begin}{\crefalias{theorem}{example}}
\AddToHook{env/remark/begin}{\crefalias{theorem}{remark}}
\AddToHook{env/question/begin}{\crefalias{theorem}{question}}
\AddToHook{env/conjecture/begin}{\crefalias{theorem}{conjecture}}

\crefname{theorem}{Theorem}{Theorems}
\crefname{lemma}{Lemma}{Lemmas}
\crefname{proposition}{Proposition}{Propositions}
\crefname{definition}{Definition}{Definitions}
\crefname{corollary}{Corollary}{Corollaries}
\crefname{assumption}{Assumption}{Assumptions}
\crefname{remark}{Remark}{Remarks}
\crefname{conjecture}{Conjecture}{Conjectures}

\newcommand{\boundary}{\partial}

\usepackage{marginnote}
    \DeclareFontFamily{U}{wncy}{}
    \DeclareFontShape{U}{wncy}{m}{n}{<->wncyr10}{}
    \DeclareSymbolFont{mcy}{U}{wncy}{m}{n}
    \DeclareMathSymbol{\Sha}{\mathord}{mcy}{"58} 

\newcommand{\Grapes}{{\mathcal{G}\mathsf{rapes}}}
\newcommand{\sha}{{\scriptstyle\Sha}}

\begin{document}
\title{On coshuffle comultiplication on configuration spaces}

\author{Byung Hee An}
\address{\parbox{\linewidth}{Department of Mathematics Education, Kyungpook National University, Daegu, Korea\\
Department of STEM Education, North Carolina State University, NC, United States, Visiting Scholar}}
\email{anbyhee@knu.ac.kr, ban2@ncsu.edu}


\begin{abstract}
We introduce a coshuffle comultiplication on the singular chain complex of configuration spaces, and we show that this structure endows the configuration space with the structure of a differential graded coalgebra ({\DGCA}). We then prove that the coshuffle comultiplication is compatible with the external product through a natural commutation relation. As an application, we investigate configuration spaces of graphs and the associated graph braid groups. In particular, for graphs of topological circumference at most 1, we prove that the singular chain complex of the configuration space is formal as a {\DGCA}. Moreover, we obtain a complete classification of the primitivity in the homology of configuration spaces of such graphs.
\end{abstract}

\maketitle

\tableofcontents

\section{Introduction}
Let \(X\) be a finite CW complex. For each \(k\ge 0\), the \emph{ordered} and \emph{unordered} \emph{\(k\)-configuration spaces} \(P_kX\) and \(B_kX\) of \(X\) is the set consisting of \(k\) distinguishable and indistinguishable distinct points in \(X\), respectively.
\begin{align*}
P_kX&=\{(x_1,\dots,x_k)\in X^k\mid x_i\neq x_j\forall i\neq j\},\\
B_kX&=\{\{x_1,\dots,x_k\}\subset X\mid x_i\neq x_j\forall i\neq j\}\cong P_kX/\mathbb{S}_k,
\end{align*}
where \(\mathbb{S}_k\) is the symmetric group on \(k\) letters acting on \(P_k X\) by permuting coordinates. The total configuration space \(BX\) is defined to be the union of all \(k\)-configuration spaces.
Note that \(B_0X\) consists of the single point \(\varnothing\) denoting the empty configuration on \(X\).

Configuration spaces have played a central role in algebraic topology since the foundational work of Fadell and Neuwirth. When \(X\) is a graph, the fundamental group of \(B_k\graf\) is known as the graph braid group, introduced around the early 2000s by Abrams and Ghrist \cite{Abrams2000,Ghrist2001}, and has since become a key object in the study of motion planning, robotics, and geometric group theory.

A distinctive feature of configuration spaces on graphs is the presence of edge stabilization maps, formalized in \cite{Ramos2018,ADCK2020}. These stabilization operations allow one to construct a finitely generated chain model \(S_\bullet(\graf)\) for \(B\graf\), known as the \'Swi\k atkowski complex, which is naturally a module over the polynomial ring generated by the edges of \(\graf\). 
In the same work, the authors established a complete classification of when this complex is formal as an edge–ring module—that is, when it is connected to its homology by a zigzag of quasi–isomorphisms of \(\ring[E]\)-modules.

On the other hand, any cochain complex of spaces carries a natural multiplication via the cup product, and dually, one may regard the chain complex as carrying a corresponding comultiplication. For configuration spaces, this duality becomes particularly meaningful because a configuration consists of many points in \(X\) that may be separated into disjoint subconfigurations. Equivalently, in the ordered context, this corresponds to taking unions of disjoint subsets of \(X\) in a shuffled manner. Motivated by this viewpoint, we introduce a coshuffle comultiplication on the singular chain complex of a configuration space.

\subsection{The results of the paper}
The central goal of this paper is to formalize this comultiplicative structure and study its interaction with the topology of configuration spaces, especially in the context of graphs. Our first main result establishes that the singular chain complex carries a counital, coaugmented, cocommutative, coassociative coalgebra({\DGCA}) structure.
\begin{theorem}[\Cref{theorem:coshuffle}]
Let \(e:C_\bullet(BX)\to C_0(B_0X)\cong\ring\) be the projection map onto the summand of degree \(0\) and weight \(0\), and let \(u:\ring\cong C_0(B_0X)\to C_\bullet(BX)\) be the map induced by the inclusion.
Then \((C_\bullet(BX),\Sha^*,e,u)\) is a {\DGCA}.
\end{theorem}

\begin{theorem}[\Cref{theorem:shuffle product on cohomology}]
Let \(X\) be a finite simplicial complex. 
If the homology group \(H_\bullet(BX)\) over \(\ring\) is torsion-free, then 
\(H_\bullet(BX)\) admits a {\DGCA} structure \((H_\bullet(BX),\Sha^*_H,e_H,u_H)\).
\end{theorem}

For a graph \(\graf\), the \emph{circumference} of \(\graf\) is the maximum of the edge lengths of embedded cycles, and the \emph{topological circumference} of \(\graf\) is the minimal circumference among all graphs \(\graf'\) homeomorphic to \(\graf\).
Then a graph with topological circumference at most \(1\) is obtained from a tree \(\sfT\) by attaching self-loops at some vertices of \(\sfT\), which can be characterized by a function \(\loops:V(\sfT)\to\bbZ_{\ge 0}\) indicating the number of loops attached to each vertex.

Due to its shape, a graph with (topological) circumference at most \(1\) is also called a \emph{bunch of grapes}, and we denote the set of bunches of grapes by \(\Grapes\).
See \Cref{figure:example of a bunch of grapes} for example.

\begin{theorem}[\Cref{theorem:formality}]
Let \(\graf=(V,E)\cong(\sfT,\ell)\) be a bunch of grapes.
Then the {\DGCA} \(S_\bullet(\graf)\) is formal, and so is \(C_\bullet(\graf)\).
\end{theorem}
That is, there are zig-zags of quasi-isomorphisms between three {\DGCA}'s,  \((C_\bullet(B\graf),\Sha^*,e,u)\), \((S_\bullet(\graf),\Sha^*_S,e_S,u_S)\), and \((H_\bullet(B\graf),\Sha^*_H,e_H,u_H)\).

Recall that for a coalgebra \((A,\Delta)\), an element \(x\in A\) is called \emph{primitive} if \(\Delta(x)=x\otimes 1+1\otimes x\).
Let \(\ring_0[E]\) be the subring of \(\ring[E]\) generated by elements of the form \(\sfe_i-\sfe_j\) over \(\ring\). See \Cref{definition:R0}.

\begin{theorem}[\Cref{theorem:primitivity}]
Let \(\graf=(V,E)\cong(\sfT,\ell)\) be a bunch of grapes.
Then \(\alpha\in H_\bullet(B\graf)\) is primitive with respect to the comultiplication \(\Sha^*_H\) if and only if \(\alpha\) is a linear combination over \(\ring\) of the following: for a fixed edge \(\sfe\in E\),
\begin{enumerate}
\item \(\sfe\cdot\varnothing\in H_0(B_0\graf)\),
\item \(p\alpha_{1jk}^\sfv\in \SL(\graf)\), or
\item \(q\beta_r^\sfv\in \SL(\graf)\),
\end{enumerate}
where \(p,q\in \ring_0[E]\) and \(\sfv\in V^\ess\).
\end{theorem}

\subsection{Outline of the paper}
The structure of the paper is as follows.

In Section~2, we introduce the notion of intertwining configuration spaces and use it to define a coshuffle comultiplication on the singular chain complex of a configuration space. This construction equips the chain complex with a natural differential graded coalgebra ({\DGCA}) structure.

In Section~3, we study the relationship between this coshuffle comultiplication and the edge stabilization maps that arise in the context of graph configuration spaces. We examine the compatibility between these two structures and clarify how the coshuffle comultiplication behaves under stabilization.

In Section~4, we apply the above framework to the class of graphs known as \emph{bunches of grapes} (graphs of topological circumference at most \(1\)). We prove that for such graphs the singular chain complex of the configuration space is formal as a {\DGCA}, and we further obtain a complete classification of primitive elements in their associated graph braid groups.

\subsection{Acknowledgments}
The author gratefully acknowledges the hospitality and support of the Department of STEM Education at North Carolina State University, where this work was carried out during the author’s visit. 
This paper was supported by Samsung Science and Technology Foundation under Project Number SSTF-BA2202-03.

\section{(Co)Shuffle (Co)multiplications}
In this section, we first consider (co)shuffle (co)multiplications on configuration spaces \(BX\).
Since every point \(\bfx\) in \(BX\) consists of number of points, by considering all possible separations of \(\bfx\) into two subsets, one can define a comultiplication \(C_\bullet(BX)\to C_\bullet(BX)\otimes C_\bullet(BX)\) on the singular chain complexes.

We also consider how the comultiplication is related with the external product of \(C_\bullet(BX)\) and the primitivity of elements in \(C_\bullet(BX)\) with respect to the defined comultiplication.

Here are some notations used throughout this paper.
\begin{enumerate}
\item A topological space \(X\) is assumed to admit a finite cell structure. 
\item A graph \(\graf=(V,E)\) is one-dimensional cell complex, whose sets of vertices and edges are denoted by \(V\) and \(E\), respectively.
\item For a sequence \(\bfk=(k_1,\dots,k_\ell)\) of nonnegative integers of length \(\ell\), the \emph{weight} of \(\bfk\) is the sum of all entries of \(\bfk\) denoted by \(\|\bfk\|\).
\item For each \(k\ge0\), we denote the set \(\{1,\dots,k\}\) by \([k]\). In particular, \([0]=\varnothing\).
\item We fix a commutative ring \(\ring\) with the unity \(1\) and a field \(\field\).
\item Differential graded coalgebra(\DGCA) is a chain complex over \(\ring\) equipped with a comultiplication, which is a counital, coaugmented, cocommutative and coassociative coalgebra.
\item Differential graded algebra(\DGA) is a cochain complex over \(\ring\) equipped with a multiplication, which is a unital, augmented, commutative and associative algebra.
\end{enumerate}

\subsection{(Co)Shuffle (co)multiplication on configuration spaces}
We introduce a configuration space consisting of configurations, which are partially ordered.

\begin{definition}[Intertwining configuration spaces]
Let \(X\) be a finite simplicial complex, and let \(\bfk=(k_1,\dots, k_\ell)\) be a sequence of nonnegative integers of length \(\ell\ge 2\).
The \emph{\(\bfk\)-intertwining configuration space} \(B_\bfk X\) over \(X\) consists of \(\ell\)-tuples \((\bfx_1,\bfx_2,\dots,\bfx_\ell)\) of configurations \(\bfx_j\in B_{k_j}X\) for all \(1\le j\le \ell\) such that \(\bfx_j\cap\bfx_{j'}=\varnothing\subset X\) for any \(j\neq j'\).
\end{definition}

Then there are canonical maps between configuration spaces as follows:
\[
\begin{tikzcd}[row sep=0]
B_{\|\bfk\|}X \ar[from=r,"\amalg_{\bfk}"'] & B_{\bfk}X \ar[r,hookrightarrow,"\iota_\bfk"] & \prod_{j=1}^\ell B_{k_j}X\\
\coprod_{j=1}^\ell\bfx_j \ar[from=r, mapsto] & (\bfx_1,\dots,\bfx_\ell) \ar[r,mapsto] & (\bfx_1,\dots,\bfx_\ell)
\end{tikzcd}
\]
Notice that the map \(\amalg\) is a covering map whose the deck transformation group is isomorphic to \(\bbS_{\bfk}\coloneqq\bbS_{\|\bfk\|}/\prod_j\bbS_{k_j}\).
Hence, for each \(\bfx\), the set \(\amalg_{\bfk}^{-1}(\bfx)\) has \(\left|\bbS_{\bfk}\right|=\frac{\|\bfk\|!}{k_1!\cdots k_\ell!}\) lifts in \(B_{\bfk}X\), where each lift is determined by an ordered decomposition of \(\bfx\) into \(\ell\) subsets \(\bfx_1,\dots,\bfx_\ell\) with \(|\bfx_j|=k_j\).

For each \(i\)-simplex \(\sigma:\triangle_i\to B_{\|\bfk\|}X\), we consider the pullback of the covering map \(B_{\bfk}X\to B_{\|\bfk\|}X\) over \(\triangle_i\), which is isomorphic to the product \(\bbS_{\bfk}\times\triangle_i\).
Then the induced map between covering spaces will be denoted by \(\amalg_{\bfk}^*\sigma:\bbS_{\bfk}\times\triangle_i\to B_{\bfk}X\), which fits into the following pullback diagram:
\[
\begin{tikzcd}[column sep=4pc]
\bbS_{\bfk} \times \triangle_i \ar[r,"\amalg_{\bfk}^*\sigma"]\ar[d,"\pr_2"'] & B_{\bfk}X\ar[d,"\amalg_{\bfk}"]\\
\triangle_i \ar[r,"\sigma"] & B_{\|\bfk\|}X
\end{tikzcd}
\]
More precisely, for each \(\bft\in\triangle_i\), the set \(\amalg_{\bfk}^*\sigma\left(\pr_2^{-1}(\bft)\right)\) consists of all lifts of \(\sigma(\bft)\).
If we fix a lift \(\tilde\sigma:\triangle_i\to B_{\bfk}X\) of \(\sigma\), then for each \(\bfs\in\bbS_\bfk\) and \(\bft\in\triangle_i\), we may denote \(\amalg_{\bfk}^*\sigma(\bfs,\bft) = \bfs\cdot\tilde\sigma(\bft)\), where \(\bfs\) acts pointwise as a deck transformation on \(B_{\bfk}X\).

For each sequence \(\bfi=(i_1,\dots,i_\ell)\) of nonnegative integers of length \(\ell\), 
we define a map between chain complexes over \(\ring\)
\[
\Sha^*_{\bfi,\bfk}:C_{\|\bfi\|}(B_{\|\bfk\|}X) 
\to \bigotimes_{j=1}^\ell C_{i_j}(B_{k_j}X)
\]
as the composition \(AW_\bfi\circ (\iota_\bfk)_*\circ \amalg_{\bfk}^*\) so that for a lift \(\tilde\sigma\) of \(\sigma\),
\[
\Sha^*_{\bfi,\bfk}(\sigma) = AW_\bfi ((\iota_\bfk)_*(\amalg_{\bfk}^*\sigma)) = 
\sum_{\bfs\in \bbS_{\bfk}} AW_\bfi((\iota_\bfk)_*(\bfs\cdot\tilde\sigma)).
\]
Here, \(AW_\bfi:C_{\|\bfi\|}(\prod_{j=1}^\ell B_{k_j}X)\to \bigotimes_{j=1}^\ell C_{i_j}(B_{k_j}X)\) is the Alexander-Whitney map, and so for the set \(\bar\cT=\{v_0,\dots,v_i\}\) of vertices of the simplex \(\triangle_i\), 
\[
AW_\bfi((\iota_\bfk)_*(\bfs\cdot\sigma)) =
\bigotimes_{j=1}^\ell \pr_j(\bfs\cdot\tilde\sigma|_{\triangle_{\bar\cT_j}}),
\]
where \(\triangle_{\bar\cT_j}\) is the face of \(\triangle_i\) whose set of vertices is \(\bar\cT_j=\{v_{i_1+\cdots+i_{j-1}}, \dots, v_{i_1+\cdots+i_j}\}\), and \(\pr_j:B_{\bfk}X\to B_{k_j}X\) is the projection.

\begin{example}
Let \(\bfi=(1,0)\), \(k=2\), and \((\sigma:\triangle_1=\{v_0,v_1\}\to B_2 I)\in C_1(B_2 I)\) defined as follows: by identifying \(0=v_0\), \(1=v_1\) and \(\triangle_1=[0,1]\),
\[
\sigma(t) = \{x_1(t), x_2(t)\}
\]
for two continuous functions \(x_1(t), x_2(t):[0,1]\to I=[0,1]\).

For \(\bfk=(2,0)\), we have a homeomorphism \(B_{\bfk}X\cong B_2X\times B_0X\) and \(\bbS_\bfk\) is trivial. 
Therefore there is a unique lift \(\tilde\sigma\) of \(\sigma\) to \(B_\bfk X\), and so
\begin{align*}
\Sha^*_{\bfi,\bfk}&=\pr_1\tilde\sigma|_{\{v_0,v_1\}}\otimes \pr_2\tilde\sigma|_{\{v_1\}}
=\sigma\otimes \varnothing.
\end{align*}

If \(\bfk=(0,2)\), then \(B_\bfk X\cong B_0X\times B_2X\) and
\begin{align*}
\Sha^*_{\bfi,\bfk}&=\pr_1\tilde\sigma|_{\{v_0,v_1\}}\otimes\pr_2\tilde\sigma|_{\{v_1\}}
=\varnothing_1\otimes \sigma(1).
\end{align*}
Here \(\varnothing_1\in C_1(B_0X)\) denotes the singular \(1\)-chain \([0,1]\to B_0X\cong\{\varnothing\}\).

For \(\bfk=(1,1)\), since \(\bbS_\bfk\) is a group of order \(2\), we have two lifts \(\tilde\sigma_1\) and \(\tilde\sigma_2\) of \(\sigma\) in \(B_\bfk X\) such that
\begin{align*}
\tilde\sigma_1&=(x_1(t),x_2(t)),&
\tilde\sigma_2&=(x_2(t),x_1(t)),
\end{align*}
and so 
\begin{align*}
\pr_1\tilde\sigma_1&=x_1(t),&
\pr_2\tilde\sigma_1&=x_2(t),&
\pr_1\tilde\sigma_2&=x_2(t),&
\pr_2\tilde\sigma_2&=x_1(t).
\end{align*}

Hence by definition of \(\Sha^*_{\bfi,\bfk}\),
\begin{align*}
\Sha^*_{\bfi,\bfk}(\sigma) &= \sum_{i=1}^2\left(
\pr_1\tilde\sigma_i|_{\{v_0,v_1\}}\otimes\pr_2\tilde\sigma_i|_{\{v_1\}}
\right)\\
&=x_1(t)\otimes x_2(1) + x_2(t)\otimes x_1(1).
\end{align*}
\end{example}

Here is another description of \(\Sha^*_{\bfi,\bfk}\) as follows:
for each \(\sigma\in C_i(B_kX)\), let us regard \(\sigma\) as a \(k\)-valued function on \(X\), or equivalently, a level-preserving embedding 
\(
\sigma:\triangle_i\times [k] \to \triangle_i\times X,
\)
such that for each \(\bft\in\triangle_i\), the image \(\sigma(\bft,[k])\) has the cardinality \(k\).
That is, the image of \(\sigma\) is homeomorphic to a disjoint union of \(k\)-copies \(\triangle_{i,1},\dots,\triangle_{i,k}\) of \(\triangle_i\) in \(\triangle_i\times X\).

Now we pick a base point \(\bft_0\in\triangle_i\) and assume that \(\sigma(\bft_0,j)\in \triangle_{i,j}\)\footnote{This process corresponds to pick a lift of \(\tilde\sigma\).}.
For a fixed \(\bfs\in\bbS_{\bfk}\), let us denote  \(\bfs(\{k_1+\dots+k_{j-1}+1,\dots,k_1+\dots+k_j\})\) by \(\bfs_j\).
Then the \(j\)-th factor of \(AW_{\bfi}((\iota_\bfk)_*(\bfs\cdot\sigma))\) is the composition 
\[
\begin{tikzcd}
\triangle_{\tilde\cT_j}\times[k_j]\ar[r,"\cong"] &
\triangle_{\tilde\cT_j}\times \bfs_j
\ar[r,hookrightarrow] &
\triangle_i\times\bfs_j \ar[r, "\sigma|_{\triangle_i\times\bfs_j}"] &
\triangle_i\times X.
\end{tikzcd}
\]
See \Cref{figure:coshuffle} for an example for \(\sigma=\left(
\begin{tikzpicture}[scale=0.5,baseline=-.5ex, yshift=.5cm]
\draw[thick] (-1,0) -- (0,-1) (0,0) -- (1,-1);
\draw[line width=5, white] (1,0) -- (-1,-1);
\draw[thick] (1,0) -- (-1,-1);
\draw[fill] (-1,0) circle (2pt) (0,0) circle (2pt) (1,0) circle (2pt);
\draw[fill] (-1,-1) circle (2pt) (0,-1) circle (2pt) (1,-1) circle (2pt);
\end{tikzpicture}\right)\).

\begin{figure}[ht]
\centering
\begin{align*}
\Sha^*_{(1,0),(3,0)}\left(\sigma\right)&=
\left(\begin{tikzpicture}[scale=0.5,baseline=-.5ex, yshift=.5cm]
\draw[thick] (-1,0) -- (0,-1) (0,0) -- (1,-1);
\draw[line width=5, white] (1,0) -- (-1,-1);
\draw[thick] (1,0) -- (-1,-1);
\draw[fill] (-1,0) circle (2pt) (0,0) circle (2pt) (1,0) circle (2pt);
\draw[fill] (-1,-1) circle (2pt) (0,-1) circle (2pt) (1,-1) circle (2pt);
\end{tikzpicture}\right)
\otimes \varnothing
\\
\Sha^*_{(1,0),(2,1)}\left(\sigma\right)
&=
\left(\begin{tikzpicture}[scale=0.5,baseline=-.5ex, yshift=.5cm]
\draw[thick] (-1,0) -- (0,-1) (0,0) -- (1,-1);
\draw[fill, gray] (1,0) circle (2pt) (-1,-1) circle (2pt);
\draw[fill] (-1,0) circle (2pt) (0,0) circle (2pt);
\draw[fill] (0,-1) circle (2pt) (1,-1) circle (2pt);
\end{tikzpicture}\right)\otimes
\left(\begin{tikzpicture}[scale=0.5,baseline=-.5ex]
\draw[fill,gray] (1,0) circle (2pt) (0,0) circle (2pt);
\draw[fill] (-1,0) circle (2pt);
\end{tikzpicture}\right)
+
\left(\begin{tikzpicture}[scale=0.5,baseline=-.5ex, yshift=.5cm]
\draw[thick] (-1,0) -- (0,-1);
\draw[line width=5, white] (1,0) -- (-1,-1);
\draw[thick] (1,0) -- (-1,-1);
\draw[fill, gray] (0,0) circle (2pt) (1,-1) circle (2pt);
\draw[fill] (-1,0) circle (2pt) (1,0) circle (2pt);
\draw[fill] (-1,-1) circle (2pt) (0,-1) circle (2pt);
\end{tikzpicture}\right)\otimes
\left(\begin{tikzpicture}[scale=0.5,baseline=-.5ex]
\draw[fill,gray] (-1,0) circle (2pt) (0,0) circle (2pt);
\draw[fill] (1,0) circle (2pt);
\end{tikzpicture}\right)
+
\left(\begin{tikzpicture}[scale=0.5,baseline=-.5ex, yshift=.5cm]
\draw[thick] (0,0) -- (1,-1);
\draw[line width=5, white] (1,0) -- (-1,-1);
\draw[thick] (1,0) -- (-1,-1);
\draw[fill, gray] (-1,0) circle (2pt) (0,-1) circle (2pt);
\draw[fill] (0,0) circle (2pt) (1,0) circle (2pt);
\draw[fill] (-1,-1) circle (2pt) (1,-1) circle (2pt);
\end{tikzpicture}\right)\otimes
\left(\begin{tikzpicture}[scale=0.5,baseline=-.5ex]
\draw[fill,gray] (-1,0) circle (2pt) (1,0) circle (2pt);
\draw[fill] (0,0) circle (2pt);
\end{tikzpicture}\right)\\
\Sha^*_{(1,0),(1,2)}\left(\sigma\right)
&=
\left(\begin{tikzpicture}[scale=0.5,baseline=-.5ex, yshift=.5cm]
\draw[thick] (-1,0) -- (0,-1);
\draw[fill, gray] (0,0) circle (2pt) (1,0) circle (2pt) (-1,-1) circle (2pt) (1,-1) circle (2pt);
\draw[fill] (-1,0) circle (2pt) (0,-1) circle (2pt);
\end{tikzpicture}\right)\otimes
\left(\begin{tikzpicture}[scale=0.5,baseline=-.5ex]
\draw[fill,gray] (0,0) circle (2pt);
\draw[fill] (-1,0) circle (2pt) (1,0) circle (2pt);
\end{tikzpicture}\right)
+
\left(\begin{tikzpicture}[scale=0.5,baseline=-.5ex, yshift=.5cm]
\draw[thick] (0,0) -- (1,-1);
\draw[fill, gray] (-1,0) circle (2pt) (1,0) circle (2pt) (-1,-1) circle (2pt) (0,-1) circle (2pt);
\draw[fill] (0,0) circle (2pt) (1,-1) circle (2pt);
\end{tikzpicture}\right)\otimes
\left(\begin{tikzpicture}[scale=0.5,baseline=-.5ex]
\draw[fill,gray] (1,0) circle (2pt);
\draw[fill] (-1,0) circle (2pt) (0,0) circle (2pt);
\end{tikzpicture}\right)
+
\left(\begin{tikzpicture}[scale=0.5,baseline=-.5ex, yshift=.5cm]
\draw[thick] (1,0) -- (-1,-1);
\draw[fill, gray] (-1,0) circle (2pt) (0,0) circle (2pt) (0,-1) circle (2pt) (1,-1) circle (2pt);
\draw[fill] (1,0) circle (2pt) (-1,-1) circle (2pt);
\end{tikzpicture}\right)\otimes
\left(\begin{tikzpicture}[scale=0.5,baseline=-.5ex]
\draw[fill,gray] (-1,0) circle (2pt);
\draw[fill] (0,0) circle (2pt) (1,0) circle (2pt);
\end{tikzpicture}\right)\\
\Sha^*_{(1,0),(0,3)}\left(\sigma\right)&=
\varnothing_1\otimes 
\left(\begin{tikzpicture}[scale=0.5,baseline=-.5ex]
\draw[fill] (-1,0) circle (2pt) (0,0) circle (2pt) (1,0) circle (2pt);
\end{tikzpicture}
\right)
\end{align*}
\caption{Example of a \(\Sha^*_{(1,0),\bfk}\).}
\label{figure:coshuffle}
\end{figure}

We also define the following maps:
\begin{align*}
\Sha^*_{i,\bfk}&=\bigoplus_{\substack{\bfi\in\bbZ_{\ge 0}^\ell\\\|\bfi\|=i}}
\Sha^*_{\bfi,\bfk}:C_i(B_{\|\bfk\|}X)\to \bigoplus_{\substack{\bfi\in\bbZ_{\ge 0}^\ell\\\|\bfi\|=i}}\bigotimes_{j=1}^\ell C_{i_j}(B_{k_j}X)
\\
\Sha^*_{\bfk}&=\bigoplus_{i\ge 0}\Sha^*_{i,\bfk}:C_\bullet(B_{\|\bfk\|}X)\to\bigotimes_{j=1}^\ell C_\bullet(B_{k_j}X),\\
\Sha^*_{\ell,i,k}&=\bigoplus_{\substack{\bfk\in\bbZ_{\ge 0}^\ell\\\|\bfk\|=k}}\Sha^*_{i,\bfk}:
C_i(B_kX)\to \bigoplus_{\substack{\bfk\in\bbZ_{\ge 0}^\ell\\\|\bfk\|=k}}\bigoplus_{\substack{\bfi\in\bbZ_{\ge 0}^\ell\\\|\bfi\|=i}} \bigotimes_{j=1}^\ell C_{i_j}(B_{k_j}X),
\\
\Sha^*_{\ell,i}&\coloneqq \bigoplus_{k\ge0} \Sha^*_{\ell,i,(k)}:C_i(BX)\to \bigoplus_{\substack{\bfi\in\bbZ_{\ge 0}^\ell\\\|\bfi\|=i}} \bigotimes_{j=1}^\ell C_{i_j}(BX),
\end{align*}
and
\begin{align*}
\Sha^*_\ell&\coloneqq \bigoplus_{i\ge 0}\Sha^*_{\ell,i}:C_\bullet(BX) \to C_\bullet(BX)^{\otimes\ell}.
\end{align*}

\begin{definition}[Coshuffle maps]
The map \(\Sha^*_{\ell}\) defined above is called the \(\ell\)\emph{-th coshuffle map}.
In particular, if \(\ell=2\), then we call \(\Sha^*_2\) a \emph{coshuffle comultiplication} and denote it simply by \(\Sha^*\).
\end{definition}

Let \(\boundary:C_\bullet(BX)\to C_\bullet(BX)\) be the boundary map. Then the boundary map \(\bar\boundary\) on \(C_\bullet(BX)^{\otimes \ell}\) is given as \(\sum_{j=1}^\ell 1^{\otimes j-1}\otimes \boundary \otimes 1^{\otimes \ell-j}\) so that for \(\sigma_j\in C_{i_j}(B_{k_j}X)\) and \(1\le j\le \ell\),
\[
\bar\boundary(\sigma_1\otimes\cdots\otimes\sigma_\ell)
=\sum_{j=1}^\ell (-1)^{i_1+\cdots+i_{j-1}}(\sigma_1\otimes\cdots\otimes\sigma_{j-1}\otimes \boundary\sigma_{i_j}\otimes\sigma_{i_{j+1}}\otimes\cdots\otimes\sigma_{\ell}).
\]

\begin{lemma}\label{lemma:chain map}
For each \(\ell\ge 2, i\ge 0, k\ge 0\) and \(\bfk\in\bbZ_{\ge 0}^\ell\), the maps \(\Sha^*_{i,\bfk}\), \(\Sha^*_{\ell,i,k}\), \(\Sha^*_{\ell,i}\), and \(\Sha^*_\ell\) are weight-preserving chain maps.
\end{lemma}
\begin{proof}
By definition, it is enough to show that \(\amalg_{\bfk}^*\) is a chain map.

Notice that \(\amalg\) is equivariant under taking subsimplices, that is, for any subsimplex \(\triangle'\subset\triangle_i\),
we have \(\amalg_{\bfk}^*\sigma|_{\bbS_{\bfk}\times \triangle'} = \amalg_{\bfk}^*(\sigma|_{\triangle'})\) and therefore \(\amalg_{\bfk}\) induces a chain map \(\amalg_{\bfk}^*:C_i(B_{\|\bfk\|}X)\to C_i(B_{\bfk}X)\).
\end{proof}

\begin{theorem}\label{theorem:coshuffle}
Let \(e:C_\bullet(BX)\to C_0(B_0X)\cong\ring\) be the projection map onto the summand of degree \(0\) and weight \(0\), and let \(u:\ring\cong C_0(B_0X)\to C_\bullet(BX)\) be the map induced by the inclusion.
Then \((C_\bullet(BX),\Sha^*,e,u)\) is a {\DGCA}.
\end{theorem}
\begin{proof}
Let \(\triangle_i\) be the \(i\)-simplex and \(\bft_0\in\triangle_i\) be a fixed point.
Suppose that \(k=k_1+k_2\) and \(\sigma\in C_i(B_kX)\). As observed earlier, summands of \(\amalg_{(k_1,k_2)}^*(\sigma)\in C_i(B_{k_1}X\times B_{k_2}X)\) are indexed by \(k_1\)-subsets of \(\amalg^{-1}(\sigma(\bft_0))\), or equivalently, \(k_2\)-subsets by considering complements.

\noindent\textbf{Cocommutativity.} 
Let \(T:B_{k_1}X\times B_{k_2}X\to B_{k_2}X\times B_{k_1}X\) be the map defined by \((\bfx_1,\bfx_2)\mapsto(\bfx_2,\bfx_1)\), which induces isomorphisms making the following diagram commutative:
\[
\begin{tikzcd}
& C_i(B_{k_1+k_2}X) \ar[rd,"\amalg_{k_2,k_1}^*"] \ar[ld,"\amalg_{k_1,k_2}^*"'] \\ C_i(B_{k_1,k_2}X) \ar[rr, leftrightarrow, "T_*"', "\cong"]\ar[d,"\iota_*"'] &
& C_i(B_{k_2,k_1}X)\ar[d,"\iota_*"]\\
C_i(B_{k_1}X\times B_{k_2}X) \ar[rr, leftrightarrow, "T_*"', "\cong"]\ar[d,"AW"'] &
& C_i(B_{k_2}X\times B_{k_1}X)\ar[d,"AW"]\\
\displaystyle\bigoplus_{\substack{i_1+i_2=i\\i_1,i_2\ge 0}} C_{i_1}(B_{k_1}X)\otimes C_{i_2}(B_{k_2}X) \ar[rr, leftrightarrow, "T_*"', "\cong"] & &
\displaystyle\bigoplus_{\substack{i_1+i_2=i\\i_1,i_2\ge 0}} C_{i_2}(B_{k_2}X)\otimes C_{i_1}(B_{k_1}X)
\end{tikzcd}
\]
Note that the upper triangle follows from the observation at the beginning of the proof.

Hence by summing up both sides in the bottom row, we have \(\Sha^* = \Sha^* T_*\), which proves the cocommutativity of \(\Sha^*\). 

\noindent\textbf{Coassociativity.} As before, summands of \((\Sha^*\otimes 1)\Sha^*(\sigma)\) in \((C_{i_1}(B_{k_1}X)\otimes C_{i_2}(B_{k_2}X))\otimes C_{i_3}(B_{k_3}X)\) are indexed by \(k_1\)-subsets of \((k_1+k_2)\)-subsets of \(\amalg^{-1}(\sigma(\bft_0))\), which are the same as the ways of decomposition of \(\amalg^{-1}(\sigma(\bft_0))\) into three ordered subsets.
In other words, 
\[
(\Sha^*\otimes 1)\Sha^* = \Sha_3^*.
\]
Similarly, we also have \((1\otimes\Sha^*)\Sha^*=\Sha_3^*\), and so \(\Sha^*\) is coassociative.

\noindent\textbf{Counitality.}
Then for each \(\sigma\in C_i(B_kX)\),
\begin{align*}
(e\otimes 1)\Sha^*(\sigma) &= 
(e\otimes 1)\Sha^*_{i,k}(\sigma)=\sigma\in
C_0(B_0X)\otimes C_i(B_kX),\\
(1\otimes e)\Sha^*(\sigma) &= 
(1\otimes e)\Sha^*_{i,k}(\sigma)=\sigma\in
C_i(B_kX)\otimes C_0(B_0X).
\end{align*}
Therefore \(e\) is a counit map.

\noindent\textbf{Coaugmentation.} It is obvious that \(u:\ring\cong C_0(B_0X)\to C_\bullet(BX)\) is injective and \(e\circ u:\ring\to\ring\) is the identity.

\noindent\textbf{Coderivation.}
Since \(\Sha^*\) is a chain map as seen in \Cref{lemma:chain map}, we have 
\[
\Sha^* \boundary = \bar\boundary\Sha^*= (\boundary\otimes 1+1\otimes\boundary)\Sha^*.
\]

For each \(\sigma\in C_i(B_kX)\), \((e\circ\boundary)(\sigma)=0\) unless \(i=1\) and \(k=0\).
However, \(C_1(B_0X)\cong \ring\) is generated by the constant map \(\sigma:[0,1]\to B_0X\cong\{*\}\), whose boundary vanishes, and therefore \((e\circ\boundary)=0\).

Finally, since \((\boundary\circ u)(1)\in C_{-1}(B_0X)=0\), we are done.
\end{proof}

\begin{remark}
Indeed, we have 
\[
(\Sha^*\otimes 1^{\ell-1})(\Sha^*\otimes 1^{\ell-2})\cdots(\Sha^*\otimes 1)\Sha^* = \Sha^*_\ell.
\]
\end{remark}

As a corollary, we have the following observation.
\begin{corollary}
The assignment \(X\mapsto (C_\bullet(BX),\Sha^*,e,u)\) is \emph{functorial}. That is, for any embedding \(X\to Y\), there is a morphism \(C_\bullet(BX)\to C_\bullet(BY)\) between {\DGCA}'s.
\end{corollary}

One another corollary is the theorem below about the induced {\DGCA} structure on the homology group.
\begin{theorem}\label{theorem:shuffle product on cohomology}
Let \(X\) be a finite simplicial complex. 
If the homology group \(H_\bullet(BX)\) over \(\ring\) is torsion-free, then 
\(H_\bullet(BX)\) admits a {\DGCA} structure \((H_\bullet(BX),\Sha^*_H,e_H,u_H)\).
\end{theorem}
\begin{proof}
By the assumption, we have an isomorphism \(\Phi:H_\bullet(BX)\otimes H_\bullet(BX)\to H(C_\bullet(BX)\otimes C_\bullet(BX))\) coming from the K\"unneth formula. Hence the comultiplication \(\Sha^*_H\) on \(H_\bullet(BX)\) is defined as the composition
\[
\begin{tikzcd}[column sep=3pc]
\Sha^*_H:H_\bullet(BX)\ar[r,"H_\bullet\Sha^*"] & H(C_\bullet(BX)\otimes C_\bullet(BX)) \ar[r,"\Phi^{-1}"] & H_\bullet(BX)\otimes H_\bullet(BX),
\end{tikzcd}
\]
which makes the conditions hold together with the counit \(e_H:H_\bullet(BX)\to H_0(B_0X)\cong\ring\) and coaugmentation \(u_H:\ring\cong H_0(B_0X)\to H_\bullet(BX)\) induced by \(e\) and \(u\) on homology groups.
\end{proof}

Now let us consider the cochain complex \((C^\bullet(BX),d)\) of configuration spaces.
Recall the chain maps \(\iota_*:C_\bullet(B_{k_1,k_2}X)\to C_\bullet(B_{k_1}X\times B_{k_2}X)\) and \(\amalg_{k_1,k_2}^*:C_\bullet(B_{k_1+k_2}X)\to C_\bullet(B_{k_1;k_2}X)\).
Then the usual cross product \(\times\) on \(C^\bullet(B_kX)\) followed by the dual of \(\iota_*\) and \(\amalg_{k_1,k_2}^*\) gives us a cochain map
\[
\begin{tikzcd}
C^\bullet(B_{k_1}X)\otimes C^\bullet(B_{k_2}X)\ar[r,"\times"] & C^\bullet(B_{k_1}X\times B_{k_2}X) \ar[r,"\iota^*"] & C^\bullet(B_{k_1,k_2}X)
\ar[r,"(\amalg_{k_1,k_2})_*"] & C^\bullet(B_{k_1+k_2}X).
\end{tikzcd}
\]

By summing up over all \(k_1,k_2\ge 0\), we have a cochain map
\begin{equation}\label{eq:shuffle}
\Sha:C^\bullet(BX)\otimes C^\bullet(BX)\to C^\bullet(BX).
\end{equation}

\begin{definition}[Shuffle product on cochains]
The map \(\Sha\) defined in \Cref{eq:shuffle} is called the \emph{shuffle product} on \(C^\bullet(BX)\).
\end{definition}

Let \(\eta:C^0(B_0X)\to C^\bullet(BX)\) be the \emph{unit} induced by the inclusion. 
Note that 
\[
C^0(B_0X)\cong \hom(C_0(B_0X),\ring)\cong \hom(\ring,\ring)\cong \ring,
\]
which is a dual of the counit \(e:C_\bullet(BX)\to C_0(B_0X)\cong\ring\).
Similarly, we define an \emph{augmentation} \(\epsilon:C^\bullet(BX)\to C^0(B_0X)\cong\ring\) to be the projection onto the factor \(C^0(B_0X)\), which is a dual of the coaugmentation \(u:\ring\cong C_0(B_0X)\to C_\bullet(BX)\).

\begin{theorem}
Let \(X\) be a finite simplicial complex.
Then the tuple \((C^\bullet(BX),\Sha,\eta,\epsilon)\) is a {\DGA}, and there is a {\DGA} structure \((H^\bullet(BX),\Sha_H,\eta,\epsilon)\) on its cohomology group \(H^\bullet(BX)\).
\end{theorem}
\begin{proof}
The first statement is a dual of \Cref{theorem:coshuffle} and we omit the proof.

The product \(\Sha_H\) on \(H^\bullet(BX)\) will be given as the composition
\[
\begin{tikzcd}[column sep=3pc]
\Sha_H:H^\bullet(BX)\otimes H^\bullet(BX)\ar[r]& H(C^\bullet(BX)\otimes C^\bullet(BX))
\ar[r,"H^\bullet\Sha"] & H^\bullet(BX),
\end{tikzcd}
\]
where the left homomorphism is a canonical map.
The unit and augmentation are induced by \(\eta\) and \(\epsilon\).
\end{proof}

\subsection{External product and primitivity}
Let \(X_1, X_2,\dots, X_\ell\) be spaces with embeddings \(\iota_j:X_j\to X\), whose images are disjoint.
Then for each \(\bfk=(k_1,\dots,k_\ell)\in\bbZ_{\ge0}^\ell\) with \(\|\bfk\|=k\), we may consider the composition \(*_\bfk\) of chain maps
\[
\begin{tikzcd}[column sep=3.5pc]
\ast_{\bfk}:\displaystyle\bigotimes_{j=1}^\ell C_\bullet(B_{k_j}X_j)\ar[r,"EZ"] &
\displaystyle C_\bullet\left(\prod_{j=1}^\ell B_{k_j}X_j\right) \ar[r, "\left(\prod\iota_j\right)_*"] & 
C_\bullet(B_{\bfk}X) \ar[r,"(\coprod_\bfk)_*"]
& C_\bullet(B_{\|\bfk\|}X).
\end{tikzcd}
\]
where \(EZ\) is the Elienberg-Zilber map.
By summing up over all sequences \(\bfk\) of length \(\ell\), we have the chain map
\[
\ast_\ell:\bigotimes_{j=1}^\ell C_\bullet(BX_j)\to C_\bullet(BX).
\]
We simply denote \(*_2\) by \(*\).
\begin{remark}
It is easy to observe that 
\[
\sigma*\tau = (-1)^{\deg\sigma\cdot\deg\tau}\tau*\sigma,
\]
and therefore one may regard the external product as the function defined on the exterior product
\(\bigwedge_{j=1}^\ell C_\bullet(BX_j)\).
\end{remark}

\begin{definition}[External products]
Let \(\sigma_j\in C_{i_j}(B_{k_j}X_j)\) for each \(1\le j\le \ell\).
The \emph{external product} \(\conv_{j=1}^\ell\sigma_j =\sigma_1*\cdots*\sigma_\ell\) via embeddings \((\iota_j)_{j=1}^\ell\) is the image of \(\sigma_1\otimes\cdots\otimes\sigma_\ell\) under the map \(\ast_\ell\)
\[
\conv_{j=1}^\ell\sigma_j = \sigma_1*\cdots*\sigma_\ell = \ast_\ell(\sigma_1\otimes\cdots\otimes\sigma_\ell).
\]

Similarly, for homology classes \(a_j\in H_{i_j}(B_{k_j}X_j)\) represented by cycles \(\sigma_j\in C_{i_j}(B_{k_j}X_j)\) for each \(1\le j\le \ell\), the external product \(a_1*\cdots*a_\ell\) is the homology class represented by the cycle \(\sigma_1*\dots*\sigma_\ell\).

We say that the external product \(c=a_1*\cdots*a_\ell\) is \emph{nontrivial} if 
\((\iota_j)_*(a_j)\)'s are linearly independent in \(H_\bullet(BX)\), and \(k_j\ge 1\) for all \(1\le j\le \ell\).
\end{definition}
\begin{remark}
When \(k_j=0\), then \(H_\bullet(B_0X_j)= H_0(B_0X_j)\cong\ring\).
\end{remark}

Then the external product is related with the comultiplication \(\Sha_\bfk^*\) as follows:
\[
\begin{tikzcd}
*_\bfk:\displaystyle\bigotimes_{j=1}^\ell C_\bullet(B_{k_j}X_j)\ar[r,"EZ"] \ar[d,"\bigotimes (\iota_j)_*"'] &
\displaystyle C_\bullet\left(\prod_{j=1}^\ell B_{k_j}X_j\right) \ar[r, "\left(\prod\iota_j\right)_*"] \ar[d,"(\prod \iota_j)_*"]& 
C_\bullet(B_{\bfk}X) \ar[r,"(\coprod_\bfk)_*"] \ar[d,equal]
& C_\bullet(B_{\|\bfk\|}X)\hphantom{:\Sha_\bfk^*}\ar[d,equal]\\
\hphantom{*_\bfk:}\displaystyle\bigotimes_{j=1}^\ell C_\bullet(B_{k_j}X) 
& \displaystyle C_\bullet\left(\prod_{j=1}^\ell B_{k_j}X\right) \ar[l, "AW"'] 
& C_\bullet(B_\bfk X) \ar[l,"\iota_*"'] 
& C_\bullet(B_{\|\bfk\|}X):\Sha_\bfk^*
\ar[l,"(\coprod_\bfk)^*"']
\end{tikzcd}
\]
Note that the left and middle squares are commutative up to homotopy, but the right square is \emph{not} commutative.

By definition of \(\Sha^*_{\bfk}\), there is a summand \(\iota_1(\sigma_1)\otimes\cdots\otimes\iota_\ell(\sigma_\ell)\) in \(\Sha^*_{\bfk}(\sigma_1*\cdots*\sigma_\ell)\) for each sequence \(\sigma_j\in C_\bullet(B_{k_j}X_j)\).
However, it is important to note that this term could be canceled in \(\Sha^*_{\bfk}(\sigma_1*\cdots*\sigma_\ell)\).

Let \(\iota_j:X_j\to X\) be embeddings as above. We define an isomorphism
\[
T:\bigotimes_{j=1}^\ell C_\bullet(BX_j)^{\otimes m}
\to 
\left(\bigotimes_{j=1}^\ell C_\bullet(BX_j)\right)^{\otimes m}
\]
as follows: let \(\sigma_{j,k}\in C_\bullet(BX_j)\) be homogeneous elements for each \(1\le j\le \ell\) and \(1\le k\le m\). Then 
\begin{align*}
&\mathrel{\hphantom{=}}T((\sigma_{1,1}\otimes\cdots\otimes\sigma_{1,m})\otimes\cdots\otimes
(\sigma_{\ell,1}\otimes\cdots\otimes\sigma_{\ell,m}))\\
&=(-1)^s(\sigma_{1,1}\otimes\cdots\otimes\sigma_{\ell,1})\otimes\cdots\otimes
(\sigma_{1,m}\otimes\cdots\otimes\sigma_{\ell,m}),
\end{align*}
where
\[
s=\sum_{j<j',k'<k} \deg\sigma_{j,k}\deg\sigma_{j',k'}.
\]

\begin{proposition}\label{proposition:commutative diagram with external product}
For each \(\ell\ge 0\), the diagram below consisting of chain maps is commutative.
\[
\begin{tikzcd}
\displaystyle\bigotimes_{j=1}^\ell C_\bullet(BX_j) \ar[rr,"*_\ell"] \ar[d,"\bigotimes \Sha_m^*"'] & & C_\bullet(BX) \ar[d,"\Sha_m^*"]\\
\displaystyle\bigotimes_{j=1}^\ell C_\bullet(BX_j)^{\otimes m} \ar[r,"T"] & \left(\displaystyle\bigotimes_{j=1}^\ell C_\bullet(BX_j)\right)^{\otimes m} \ar[r,"*_\ell^{\otimes m}"]& C_\bullet(BX)^{\otimes m}
\end{tikzcd}
\]
\end{proposition}
\begin{proof}
Let \(\sigma_j:\triangle_{i_j}\to B_{k_j}X_j\) for each \(1\le j\le \ell\) and let \(i=i_1+\cdots+i_\ell\) and \(k=k_1+\cdots+k_\ell\).
We denote the set of ordered vertices of \(\triangle_{i_j}\) by \(\bar\cV_j=(v_0^j,\dots,v_{i_j}^j)\).
Then the Eilenberg-Zilber map decomposes the domain \(\prod_{j=1}^\ell \triangle_{i_j}\) of the product \(\prod_{j=1}^\ell \sigma_j\) into \(i\)-dimensional simplices \(\triangle_\cT\) indexed by the set of shuffles \(\cT=(w_1,\dots,w_i)\) on ordered sets
\(\cV_1=(v_1^1,\dots,v_{i_1}^1), \cdots, \cV_\ell=(v_1^\ell,\dots,v_{i_\ell}^\ell)\).
Hence, we may identify the set of vertices of \(\triangle_\cT\) with the ordered set \(\bar\cT=(w_0,w_1,\dots,w_i)\), where \(w_0\) is the common vertex of all \(\triangle_\cT\)'s and corresponds to the point \((v_0^1,v_0^2,\dots,v_0^\ell)\in \prod_{j=1}^\ell\triangle_{i_j}\).

Let \(\sigma_\cT:\triangle_\cT\to B_kX\) be the restriction of \(\prod_{j=1}^\ell \sigma_j\) to the simplex \(\triangle_\cT\), and let \(\bfk'=(k_1',\dots,k_m')\) and \(\bfi'=(i_1',\dots,i_m')\) be sequences of nonnegative integers of length \(m\ge2\) with \(\|\bfk'\|=k\) and \(\|\bfi'\|=i\).
Then for each lift \(\tilde\sigma_\cT:\triangle_\cT\to B_{\bfk'}X\) of \(\sigma_\cT\), 
\[
\Sha^*_{\bfi',\bfk'}(\sigma_\cT)=\sum_{\bfs\in\bbS_{\bfk'}}AW_{\bfi'}((\iota_{\bfk'})_*)(\bfs\cdot\tilde\sigma_\cT).
\]

Define subsequences of \(\bar\cT\) as follows:
\begin{align*}
\cT_1&=(w_1,\dots, w_{i_1'}),&
\cT_2&=(w_{i_1'+1},\dots, w_{i_1'+i_2'}),& &\cdots&
\cT_m&=(w_{i_1'+\cdots+i_{m-1}'+1},\dots, w_i),\\
\bar\cT_1&=(w_0,\dots, w_{i_1'}),&
\bar\cT_2&=(w_{i_1'},\dots, w_{i_1'+i_2'}),& &\cdots&
\bar\cT_m&=(w_{i_1'+\cdots+i_{m-1}'},\dots, w_i).
\end{align*}
Then each \(\bar\cT_{j'}\) corresponds to an \(i'_{j'}\)-dimensional face \(\triangle_{\bar\cT_{j'}}\) of \(\triangle_\cT\).
Let us consider ordered subsets
\begin{align*}
\cV_{1,j'} &= \cV_1\cap \cT_{j'},&
\cV_{2,j'} &= \cV_2\cap \cT_{j'},& &\cdots&
\cV_{\ell,j'} &= \cV_\ell\cap \cT_{j'}.
\end{align*}
Obviously, \(\cT_{j'}\) could be realized as a shuffle of \(\cV_{1,j'},\dots,\cV_{\ell,j'}\).

On the other hand, since the set \(\cV_j\) of vertices in \(\triangle_{i_j}\) is decomposed into the subsets \(\cV_{j,1}, \dots, \cV_{j,m}\), there is a summand of \(\bigotimes \Sha^*_m(\bigotimes \sigma_j)\) which is of the form
\[
(\sigma_{1,1}\otimes\cdots\sigma_{1,m})\otimes\cdots\otimes(\sigma_{\ell,1}\otimes\cdots\sigma_{\ell,m}),
\]
where \(\sigma_{j,j'}:\triangle_{\bar\cV_{j,j'}}\to C_\bullet(BX_j)\). 
Then by applying the isomorphism \(T\), we need to apply the external product \(*_\ell\) on \((\sigma_{1,j'}\otimes\cdots\otimes\sigma_{\ell,j'})\), each of whose summands has the simplex indexed by shuffles on \(\cV_{1,j'},\dots, \cV_{\ell,j'}\).

The sign matches since the sign of each shuffle \(\cT\) can be decomposed into a block shuffle on \(\{\cV_{j,j'}\}\) by \(T\) and a shuffle on \(\cV_{1,j'},\dots, \cV_{\ell,j'}\).
\end{proof}

\begin{definition}[Primitive elements]
For a coalgebra \((A,\Delta)\), an element \(x\in A\) is called \emph{primitive} if \(\Delta(x)=x\otimes 1+1\otimes x\). Namely, \(x\) can not be decomposed to two nontrivial elements with respect to the comultiplication \(\Delta\).
\end{definition}

\begin{theorem}\label{theorem:nonprimitive external product}
For \(1\le j\le \ell\), let \(\{\iota_j:X_j\to X\}\) be a sequence of embeddings whose images are disjoint.
If the external product \(a=a_1*\cdots*a_\ell\) is nontrivial for primitive \(a_j\in H_{i_j}(B_{k_j}X_j)\), then it is not primitive.
\end{theorem}
\begin{proof}
Suppose that each \(a_j\in H_\bullet(BX_j)\) is represented by \(\sigma_j\in C_{i_j}(B_{k_j}X_j)\) for some \(k_j\ge 1\).
Then \(a_1*\cdots*a_\ell\) is represented by a cycle
\[
\begin{tikzcd}
\sigma:\prod_{j=1}^\ell \triangle_{i_j} \ar[r,"\prod \sigma_j"] &
\prod_{j=1}^\ell B_{k_j}X_j \ar[r,"\prod \iota_j"] &
B_\bfk X \ar[r] & B_{k} X,
\end{tikzcd}
\]
where \(\bfk=(k_1,\dots, k_\ell)\) and \(k=k_1+\dots+k_\ell\).

Let \(\bfk'=(k_1',\dots,k_\ell')\) be a sequence of positive integers with \(\|\bfk'\|=k\).
Then for any fixed element in \(\bbS_{\bfk'}\), we determine a lift \(\tilde\sigma\) of \(\sigma\) to \(B_{\bfk'}X\) and consider the following composition:
\[
\begin{tikzcd}
\tilde\sigma_{j,j'}:\triangle_{i_j}\ar[r, hookrightarrow] &\prod_{j=1}^\ell \triangle_{i_j}\ar[r,"\tilde\sigma"] & B_{\bfk'}X\ar[r,"\operatorname{Proj}_{j'}"]& B_{k_{j'}}X,
\end{tikzcd}
\]
where the left inclusion is canonical since each simplex is pointed by the very first vertex.
By definition of the \(\ell\)-th coshuffle \(\Sha^*_\ell\), each summand of \(\Sha^*_\ell(\sigma)\) is of the form 
\[
\left(
\prod_{j=1}^\ell\tilde\sigma_{j,1}
\right)\otimes
\left(
\prod_{j=1}^\ell\tilde\sigma_{j,2}
\right)\otimes\cdots\otimes
\left(
\prod_{j=1}^\ell\tilde\sigma_{j,\ell}
\right).
\]

Notice that the lift \(\tilde\sigma\) distributes the image of \(\sigma_j\) into \(\ell\) subsets, and therefore the image of \(\tilde\sigma_{j,j'}\) in \(B_{k_{j'}}X\) may have stationary points.
More precisely, let \(v_0^j\) be the very first vertex of the simplex \(\triangle_{i_j}\) for each \(1\le j\le \ell\), and let \(\bft_0=(v_0^1,\dots,v_0^\ell)\in\prod_{j=1}^\ell\triangle_{i_j}\).
If we consider the image \(\bfx_0=\{x_1^0,\dots,x_k^0\}=\sigma(\bft_0)\), then each lift \(\tilde\sigma\) gives us an ordered partition \((\bfx_1',\dots,\bfx_\ell')\) of \(\bfx_0\), where \(|\bfx'_{j'}|=k_{j'}\) for each \(1\le j'\le \ell\).

Let \(k_{j,j'}=|\bfx'_{j'}\cap X_j|\). Then the image of \(\tilde\sigma_{j,j'}\) has at least \((k'_{j'} - k_{j,j'})\) stationary points.
Hence by ignoring stationary points or considering points in \(X_j\), we have a lift \(\tilde\sigma_j\) of \(\sigma_j\) to \(B_{\bfk_j}X_j\), where \(\bfk_{j}=(k_{j,1},k_{j,2},\dots,k_{j,\ell})\).
In other words, we have the commutative diagram
\[
\begin{tikzcd}[column sep=3pc]
\triangle_{i_j} \ar[r, "\tilde\sigma_{j,j'}"] \ar[d, "\tilde\sigma_j"'] & B_{k_{j'}X} \\
B_{\bfk_{j}X_j} \ar[r, "\pr_{j'}"] & B_{k_{j,j'}X_j}, \ar[u]
\end{tikzcd}
\]
where the right arrow is the composition of the inclusion \(\iota_j:X_j\to X\) and the map taking union with \(\bfx'_{j'}\setminus X_j\).

In particular, since \(\sigma_j\) is primitive, its homology class vanishes unless all images of \(\sigma_j\) is contained in a single subset of the partition so that \(k_{j,j'}\ge k_j\).
Since \(k_{j,j'}\le k'_{j'}\) and \(\|\bfk\|=\|\bfk'\|=k\), for any \(1\le j\le \ell\), there exists \(1\le j'\le \ell\) such that \(k_j=k_{j,j'}=k'_{j'}\).
That is, two sequences \(\bfk\) and \(\bfk'\) are the same up to a permutation, say \(\tau\).
Then there is a nontrivial summand \((\iota_{\tau(1)})_*(a_{\tau(1)})\otimes\cdots\otimes (\iota_{\tau(\ell)})_*(a_{\tau(\ell)})\) in the \(\ell\)-th coshuffle \(\Sha^*_\ell(\sigma)\), and we are done.
\end{proof}

\section{Comultiplication on configuration spaces of graphs}

\subsection{The \'Swi\k atkowski complexes for graphs}
From now on, we focus on the case when \(X\) is a graph
\( \graf=(V,E) \), which is assumed to have no bivalent vertices.
For a sake of convenience, we fix a linear order on \(V\).

\begin{example}\label{example:trivial}
Let \(\graf\) be a graph having no essential vertex. Then \(\graf\) is homeomorphic to either an interval \(I\) or a circle \(S^1\), and so we may assume that \(\graf\) has only one edge \(\sfe\).
If \(\graf\cong I\), then
\[
H_\bullet(BI) \cong H_0(BI) \cong \ring[\sfe],
\]
and if \(\graf\cong S^1\), then 
\[
H_\bullet(BS^1)\cong H_0(BS^1)\oplus H_1(BS^1) \cong \ring[\sfe]\oplus \ring[\sfe]\langle\beta\rangle,
\]
where \(1\in H_0(B_0\graf)\cong\ring\) and \(\beta\in H_1(B_1S^1)\cong H_1(S^1)\) are the homology classes represented by the point \(*\cong B_0\graf\) and the circle \(S^1\cong B_1S^1\), respectively.
In particular, we have
\[
\Sha^*(1) = 1\otimes1\quad\text{ and }\quad
\Sha^*(\beta) = \beta\otimes1 + 1\otimes\beta.
\]
\end{example}

For each edge \(\sfe\), we may assume that a parametrization \(\phi_\sfe:I=[0,1]\to\sfe\subset\graf\) for \(\sfe\) is fixed, where the restriction \(\phi_\sfe|_{(0,1)}\) to the unit open interval \((0,1)\) is a homeomorphism.
A \emph{half-edge} \(\sfh\) adjacent to a vertex \(\sfv\) is (an image of) a restriction \(\phi_{\sfe}|_{[0,1/2]}\) or \(\phi_{\sfe}|_{[1/2,1]}\) for some edge \(\sfe\) such that \(\phi_\sfe(0)=\sfv\) or \(\phi_\sfe(1)=\sfv\), respectively.
In this case, we denote \(\sfe\) and \(\sfv\) by \(\sfe(\sfh)\) and \(\sfv(\sfh)\), which are the edge that \(\sfh\) is contained in and the vertex that \(\sfh\) contains, respectively.

We denote the set of all half-edges whose vertex is \(\sfv\) by \(\cH_\sfv\).
Then for each vertex \(\sfv\), the number of half-edges adjacent to \(\sfv\) is the same as the valency (or degree) of the vertex \(\sfv\).

\begin{definition}[The \'Swi\k{a}tkowski complex]
Let \(\graf=(V,E)\) be a graph. For each vertex \(\sfv\in V^\ess\), we define a finitely generated, bigraded, free \(\ring[E]\)-module \((S_\sfv,\partial_\sfv)\) as follows:
\begin{align}\label{eq:Swiatkowski}
S_\sfv &= \ring[E]\oplus \ring[E]\langle\sfv\rangle\oplus \ring[E]\langle\cH_\sfv\rangle,&
\partial_\sfv(\gamma) &= 
\begin{cases}
0 & \text{if }\gamma=\sfv\text{ or }\gamma\in\ring[E];\\
\sfe(\sfh)-\sfv & \text{if }\gamma=\sfh\in\cH_\sfv.
\end{cases}
\end{align}
The bigrading consists of the \emph{homological degree} and the \emph{braid index}, given as follows: for each \(r\in \ring\), \(\sfv\in V\),
\(\sfh\in\cH_\sfv\), and \(\sfe\in E\),
\begin{align*}
|r|&=(0,0),&
|\sfv|&=(0,1),&
|\sfh|&=(1,1),&
&\text{and}&
|\sfe|&=(0,1).
\end{align*}

The \emph{\'Swi\k{a}tkowski complex} is a finitely generated, bigraded, free \(\ring[E]\)-module \(S_\bullet(\graf)\) defined as the exterior product over \(\ring[E]\)
\[
\left(S_\bullet(\graf),\partial_S\right)=\bigwedge_{\sfv\in V^\ess} \left(S_\sfv,\partial_\sfv\right).
\]
\end{definition}

\begin{remark}
This is a partially reduced version. The original definition of \(S_\bullet(\graf)\) is the tensor product that ranges over all vertices.
\end{remark}

\begin{remark}
We may omit the notation "\(\wedge\)" to denote elements in \(S_\bullet(\graf)\).
\end{remark}

Let \(\sfe\) be an edge of \(\graf\).
For \(\bfx\in B_n\graf\), let \(\bfx_\sfe=\bfx\cap\sfe\). Then \(\bfx_\sfe=\{x_1,\dots,x_k\}\in B_k\sfe\) for some \(k\ge 0\).
By using the parametrization \(\phi_\sfe:I\to\sfe\) of \(\sfe\), we consider \(\bft=\phi_{\sfe}^{-1}(\bfx_\sfe)=\{t_1,\dots, t_\ell\}\in B_\ell I\), where \(\ell\) is either \(k\) or \(k+1\) if \(\sfe\) is a loop adjacent to a vertex \(\sfv\) in \(\graf\) and \(\sfv\in\bfx_\sfe\).
Without loss of generality, we may assume that \(0\le t_1<\dots<t_\ell\le 1\).

We define the set \(\bft'=\{t'_1,\dots,t'_{\ell+1}\}\in B_{\ell+1}[0,1]\) by \( t_i' = (t_{i-1}+t_i)/2 \) for \( 1\le i\le \ell+1 \), where we set \(t_0=0\) and \(t_{\ell+1}=1\) for convenience. 
Notice that \(t_1=0\) if and only if \(t_1'=0\), and \(t_\ell=1\) if and only if \(t'_{\ell+1}=1\). Therefore, no matter whether \(\sfe\) is a loop or not, the image of \(\bft'\) under \(\phi_\sfe\) defines a configuration \(\bfx'_\sfe=\{x_1,\dots,x_{k+1}\}\in B_{k+1}\sfe\).

The \emph{edge stabilization} \(\sfe\cdot(-):B_k\graf\to B_{k+1}\graf\) on \(\sfe\) is the map defined by 
\[
\sfe\cdot\bfx = \left(\bfx\setminus \bfx_\sfe\right) \cup \bfx_\sfe'.
\]
It is known that the edge stabilization on \(\sfe\) is a continuous map, which allows us to regard the singular chain complex \(C_\bullet(BX)\) as a module over the multivariate polynomial ring \(\ring[E]\).

One can consider a canonical inclusion \(f:S_\bullet(\graf)\to C_\bullet(B\graf)\), which is \(\ring[E]\)-linear and defined by the following rules: let \(\sfv\) be a vertex, \(\sfe\) be an edge, \(\sfh\) and \(\sfh'\) be two half-edges of \(\sfe\), and \(\alpha\in S_\bullet(\graf)\). Then
\begin{enumerate}
\item \(f(1)\coloneqq\varnothing\in C_0(B_0\graf)\),
\item \(f(\sfv)\coloneqq\sfv\in C_0(B_1\graf)\) and \(f(\sfh_0)\coloneqq\sfh_0\in C_1(B_1\graf)\),
\item \(f(\sfe\alpha)\coloneqq\sfe\cdot f(\alpha)\) and \(f(\sfv\alpha)\coloneqq \sfv\ast f(\alpha)\),
\item \(f(\sfh\sfh'):[0,1]^2\to B_2(\sfe)\subset B_2(\graf)\) is defined as 
\[
f(\sfh\sfh')(s,t)=\phi_{\sfe}(\{s(2-t)/4, 1-(2-s)t/4\})
\]
if \(\sfv(\sfh)\neq\sfv(\sfh')\), \(\sfh=\phi_{\sfe}(s/2):[0,1]\to\sfe\) and \(\sfh'=\phi_{\sfe}(1-s/2):[0,1]\to\sfe\).
\item for \(\alpha=\sfh_1\dots\sfh_k\) with \(\sfh,\sfh'\neq\sfh_i\) for all \(i\), 
\begin{enumerate}
\item \(f(\alpha\sfh)\coloneqq f(\alpha)\ast f(\sfh)\),
\item \(f(\alpha\sfh\sfh')\coloneqq f(\alpha)\ast f(\sfh\sfh')\) if \(\sfv(\sfh)\neq\sfv(\sfh')\).
\end{enumerate}
\end{enumerate}

Note that \(f(\sfh\sfh')(0,t)=\sfv(\sfh)*\sfh'\) and \(f(\sfh\sfh')(s,0)=\sfh*\sfv(\sfh')=\sfv(\sfh')\ast\sfh\), \(f(\sfh\sfh')(1,t)=\sfe\cdot\sfh'\) and \(f(\sfh\sfh')(s,1)=\sfe\cdot\sfh\).

\begin{lemma}
The map \(f\) is a chain map between \(\ring[E]\)-modules.
\end{lemma}
\begin{proof}
Since \(f\) is \(\ring[E]\)-linear and the differential \(\partial\) on \(C_\bullet(B\graf)\) is a derivation with respect to the external product, we only need to check \(\partial f(\alpha)=f(\partial\alpha)\) for \(\alpha=1, \sfv, \sfh\) and \(\sfh\sfh'\) if \(\sfv(\sfh)\neq\sfv(\sfh')\).

For \(\alpha=1,\sfv\), this is trivial, and for \(\alpha=\sfh\), \(f(\partial_S\sfh_1)=\sfe(\sfh_1)-\sfv(\sfh_1)=\partial(\sfh_1)=\partial f(\sfh_1)\), and we are done.

Let \(\sfv(\sfh)=\sfv\) and \(\sfv(\sfh')=\sfv'\) with \(\sfv\neq\sfv'\).
Then 
\begin{align*}
f(\partial_S(\sfh\sfh'))
&=f((\partial_S\sfh)\sfh'-\sfh(\partial_S\sfh'))\\
&=f((\sfe-\sfv)\sfh'-\sfh(\sfe-\sfv'))\\
&=\sfe\cdot\sfh'-\sfv\ast\sfh'-\sfe\cdot\sfh+\sfv'\ast\sfh\\
&=f(\sfh\sfh')(1,t)-f(\sfh\sfh')(0,t)-f(\sfh\sfh')(s,1)+f(\sfh\sfh')(s,0)\\
&=\partial f(\sfh\sfh').\qedhere
\end{align*}
\end{proof}

Indeed, the map \(f\) is a chain equivalence.
\begin{theorem}\cite[Theorems~2.10 and 1.1]{ADCK2020}\label{thm:ADCK2020}
Let \(\graf=(V,E)\) be a graph. Then the singular chain complex \(C_\bullet(B\graf)\) is an \(\ring[E]\)-module, chain homotopy equivalent to \(S_\bullet(\graf)\) as \(\ring[E]\)-modules.
\end{theorem}
\begin{proof}
Since both \(C_\bullet(B\graf)\) and \(S_\bullet(\graf)\) are free \(\ring[E]\)-modules, we may identify \(C_\bullet(B\graf)\) and \(S_\bullet(\graf)\) with direct sums of the kernels and images of their differentials so that
\begin{align*}
C_i(B\graf)&\cong \image(\partial_i:C_i(B\graf)\to C_i(B\graf))\oplus
\ker(\partial_i:C_i(B\graf)\to C_{i-1}(B\graf))\\
S_i(\graf)&\cong \image(\partial:S_{i+1}(\graf)\to S_i(\graf))\oplus
\ker(\partial:S_i(\graf)\to S_{i-1}(\graf))
\end{align*}
as free \(\ring[E]\)-modules. Then for each \(i\ge 0\), we have two free resolutions over \(\ring\)
\[
\begin{tikzcd}
\image(\partial_{i+1}:C_{i+1}(B\graf)\to C_i(B\graf))\ar[r]\ar[d,"\exists g_{i+1}", xshift=.5ex]&
\ker(\partial_i:C_i(B\graf)\to C_{i-1}(B\graf))\ar[r]\ar[d,"\exists g_i", xshift=.5ex]&0\\
\image(\partial:S_{i+1}(\graf)\to S_i(\graf))\ar[r]\ar[u,"f",xshift=-.5ex]&
\ker(\partial:S_i(\graf)\to S_{i-1}(\graf))\ar[r]\ar[u,"f",xshift=-.5ex]&0
\end{tikzcd}
\]
of \(H_i(B\graf)\cong H_i(S_\bullet(\graf))\), which are chain equivalent.
Hence \(g\) is defined by collecting the chain homotopy inverses \(g_i\).

Indeed, \(g\) is a chain map over \(\ring[E]\) module since \(\partial(\sfe\cdot\sigma)=\sfe\cdot\partial\sigma\) for every \(\sigma\in C_\bullet(BX)\) as desired.
\end{proof}

\begin{example}\label{example:computation of g}
Let \(\sigma:[0,1]^k\to B_n\graf\) be a singular \(k\)-chain.
If \(n=0\), then since \(B_0\graf=\{\varnothing\}\) and the summand of bigrading \((k,0)\) of \(S_\bullet(\graf)\) is \(\ring\), which is concentrated in degree \(k=0\), we have \(g(\sigma)=0\) unless \(k=0\).

Similarly, for \(n=1\), the bigrading \((k,1)\)-summand of \(S_\bullet(\graf)\) is 
\[
\ring[E]\oplus \ring\langle V^\ess\rangle\oplus \left(\bigoplus_{\sfv\in V^\ess}\ring\langle \cH_\sfv\rangle\right).
\]
Hence \(g(\sigma)=\sfh\) if \(\sigma(0)=\sfv(\sfh)\) and \(\sigma(1)\) is contained in the interior of the edge \(\sfe(\sfh)\) for some \(\sfh\). In particular, if the image of \(\sigma\) is contained in the interior of an edge \(\sfe\), then \(g(\sigma)=0\).
\end{example}

Moreover, it is easy to check that each \(g_i\) preserves the external product and so does \(g\). That is, \(g\circ *=*\circ(g\otimes g)\).
Therefore we define the comultiplication \(\Sha^*_S:S_\bullet(\graf)\to S_\bullet(\graf)\otimes S_\bullet(\graf)\) by 
\[
\Sha^*_S(\alpha) = (g\otimes g)\Sha^*(f(\alpha)).
\]

Recall the counit \(e:C_\bullet(B\graf)\to C_0(B_0\graf)\cong \ring\) and coaugmentation \(u:\ring\cong C_0(B_0\graf)\to C_\bullet(B\graf)\) for \(C_\bullet(B\graf)\).
We define the counit \(e_S\) and coaugmentation \(u_S\) for \(S_\bullet(\graf)\) in a similar way since \(S_\bullet(\graf)\) has a direct summand isomorphic to \(\ring\).

In summary, we have the following theorem.
\begin{theorem}
Two {\DGCA}'s \((C_\bullet(B\graf),\Sha^*,e,u)\) and \((S_\bullet(\graf),\Sha^*_S, e_S, u_S)\) are isomorphic.
\end{theorem}

\subsection{Comultiplication \(\Sha^*_S\) on \(S_\bullet(\graf)\)}
Now we will compute \(\Sha^*_S\) explicitly as follows:
Let us regard \(S_\bullet(\graf)\otimes S_\bullet(\graf)\) as a module over \(\ring[E]\otimes\ring[E]\). That is, for each \(\sfe\otimes\sfe'\in \ring[E]\otimes\ring[E]\) and \(\alpha\otimes\alpha'\in S_\bullet(\graf)\otimes S_\bullet(\graf)\), we have 
\[
(\sfe\otimes\sfe')\cdot(\alpha\otimes\alpha') = \sfe\cdot\alpha\otimes\sfe'\cdot\alpha'.
\] 
We define a ring homomorphism \(\sha^*:\ring[E]\to\ring[E]\otimes\ring[E]\) by \(\sha^*(\sfe) = \sfe\otimes1+1\otimes\sfe\) for each \(\sfe\in E\).
Via the map \(\sha^*\), we may regard \(S_\bullet(\graf)\otimes S_\bullet(\graf)\) as the \(\ring[E]\)-module as well. That is, for each \(\sfe\in\ring[E]\),
\[
\sfe\cdot(\alpha\otimes\alpha') \coloneqq \sha^*(\sfe)\cdot(\alpha\otimes\alpha')
=(\sfe\otimes1+1\otimes\sfe)\cdot(\alpha\otimes\alpha').
\]

\begin{proposition}
The comultiplication \(\Sha^*_S:S_\bullet(\graf)\to S_\bullet(\graf)\otimes S_\bullet(\graf)\) is the \(\ring[E]\)-linear map completely determined by the following rules:
\begin{enumerate}
\item \(\Sha_S^*(1) = 1\otimes1\),
\item \(\Sha_S^*(\sfv) = \sfv\otimes 1+1\otimes\sfv\) and \(\Sha_S^*(\sfh) = \sfh\otimes 1+1\otimes\sfh\) for each \(\sfv\in V^\ess\) and \(\sfh\in\cH_\sfv\), and
\item 
\(
\Sha_S^*(\sfh_1\cdots\sfh_k)=
\sum_{I\subset[k]}(-1)^{s(I)} \left(\conv_{i\in I}\sfh_i\right)
\otimes\left(\conv_{j\not\in I}\sfh_j\right),
\)
\end{enumerate}
where \(s(I,J)=\#\{(i,j):i\in I,j\in J, i>j\}\).
\end{proposition}
\begin{proof}
Since \(f, \Sha^*\) and \(\Sha^*_S\) are \(\ring\)-linear, and by definition of \(f\), we only need to consider \(\alpha=\sfh_1\dots\sfh_k\).

As before, we use induction on \(k\ge1\).
Let us follow the definition of \(\Sha^*\) on \(C_\bullet(B\graf)\).

If \(k=1\), then since \(f(\sfh_1)\in C_1(B_1\graf)\), \(\bfk\) is either \((1,0)\) or \((0,1)\).
\begin{align*}
\Sha_{(1,0)}^*(f(\sfh_1))&=\Sha^*(\sfh_1)=\sfh_1\otimes \pr_2(\sfe(\sfh_1)) + \sfv(\sfh_1)\otimes \pr_2(\sfh_1).
\end{align*}
Here, \(\pr_1:B_{(1,0)}(\graf)\cong B_1\graf\times B_0\graf\to B_1\graf\) and \(\pr_2:B_{(1,0)}(\graf)\to B_0\graf\) and \(B_0\graf=\{\varnothing\}\).
Since \(\pr_2(\sfe(\sfh_1))\in C_0(B_0\graf)\cong\ring\langle\varnothing\rangle\) and \(\pr_2(\sfh)\in C_1(B_0\graf)\), we have \(g(\pr_2(\sfe(\sfh_1)))=\varnothing\) and \(g(\pr_2(\sfh_1))=0\), and therefore
\[
(g\otimes g)\Sha^*_{(1,0)}(f(\sha_1))=g(\sfh_1)\otimes 1=\sfh_1\otimes1.
\]
Similarly, \((g\otimes g)\Sha^*_{(0,1)}(f(\sha_1))=1\otimes g(\sfh_1)=1\otimes \sfh_1\), and so
\begin{align*}
\Sha^*_S(\sfh_1)&=(g\otimes g)\Sha^*_{(1,0)}(f(\sha_1))+(g\otimes g)\Sha^*_{(0,1)}(f(\sha_1))\\
&= \sfh_1\otimes 1+1\otimes \sfh_1\\
&=\sum_{I\subset[1]=\{1\}}(-1)^{s(I,[1]\setminus I)}\left(\bigwedge_{i\in I}\sfh_i\right)\otimes
\left(\bigwedge_{j\in[k]\setminus I}\sfh_j\right).
\end{align*}

It is important to note that \(\Sha^*(\sfh)\) itself on \(C_\bullet(B\graf)\) is not the same as \(\sfh\otimes1+1\otimes\sfh\), but the same up to kernel of \((g\otimes g)\).
We simply denote elements in the kernel by \(\xi\).

Let \(\sfh,\sfh'\) be two half-edges of an edge \(\sfe\), and let \(\sfv(\sfh)=\sfv\) and \(\sfv(\sfh')=\sfv'\) as before.
Notice that \(f(\sfh\sfh')\in C_2(B_2\sfe)\subset C_2(B_2\graf)\).
We cut the domain \([0,1]\times[0,1]\) of \(f(\sfh\sfh')\) into two simplices 
\[
\sigma:\triangle_1=\{v_0,v_1,v_2\}\to B_2\sfe\quad\text{ and }\quad\tau:\triangle_2=\{v_0,v_3,v_2\}\to B_2\sfe,
\]
where \(v_i\in[0,1]\times[0,1]\) are vertices
\begin{align*}
v_0&=(0,0),&
v_1&=(1,0),&
v_2&=(1,1),&
v_3&=(0,1).
\end{align*}
Then \(f(\sfh\sfh')=\sigma-\tau\).

For \(\bfk=(k_1,k_2)=(2,0)\), we have 
\begin{align*}
\Sha^*_{(2,0)}(\sigma)&=\sigma\otimes\varnothing + \sigma|_{\{v_0,v_1\}}\otimes \pr_2(\sigma|_{\{v_1,v_2\}}) + \sigma|_{\{v_0\}}\otimes \pr_2(\sigma)\\
\Sha^*_{(2,0)}&=\tau\otimes\varnothing + \tau|_{\{v_0,v_3\}}\otimes \pr_2(\tau|_{\{v_3,v_2\}}) + \tau|_{\{v_0\}}\otimes \pr_2(\tau).
\end{align*}
Since \(\pr_2:B_{(2,0)}\graf\cong B_2\graf\times B_0\graf\to B_0\graf=\{\varnothing\}\), we have
\[
g(\pr_2(\sigma|_{\{v_1,v_2\}}))=g(\pr_2(\sigma))=0
\]
as seen in \Cref{example:computation of g}, and therefore 
\begin{align*}
(g\otimes g)\Sha^*_{(0,2)}(f(\sfh\sfh'))&=(g\otimes g)\Sha^*_{(2,0)}(\sigma-\tau)\\
&=g(\sigma-\tau)\otimes g(\varnothing)\\
&=g(f(\sfh\sfh'))\otimes 1\\
&=(\sfh\sfh')\otimes 1.
\end{align*}

Similarly,
\begin{align*}
(g\otimes g)\Sha^*_{(0,2)}(f(\sfh\sfh'))=(g\otimes g)\Sha^*_{(0,2)}(\sigma-\tau)
=1\otimes g(f(\sfh\sfh'))=1\otimes(\sfh\sfh').
\end{align*}

Finally, for \(\bfk=(1,1)\), there are two lifts of \(\sigma\) and \(\tau\), denoted by \(\tilde\sigma_i\) and \(\tilde\tau_i\) for \(i=1,2\). Then 
\begin{align*}
\tilde\sigma_1(s,t)&=(\phi_\sfe(s(2-t)/4),\phi_\sfe(1-(2-s)t/4)),&
\tilde\sigma_2(s,t)&=(\phi_\sfe(1-(2-s)t/4),\phi_\sfe(s(2-t)/4))
\end{align*}
on \(\triangle_1\), and two maps \(\tilde\tau_1\) and \(\tilde\tau_2\) are defined by the same way as above on \(\triangle_2\).
\begin{align*}
\Sha^*_{(1,1)}(\sigma)&=\pr_1\tilde\sigma_1\otimes\pr_2\tilde\sigma_1|_{\{v_2\}}
+\pr_1\tilde\sigma_1|_{\{v_0,v_1\}}\otimes\pr_2\tilde\sigma_1|_{\{v_1,v_2\}}
+\pr_1\tilde\sigma_1|_{\{v_0\}}\otimes\pr_2\tilde\sigma_1\\
&\mathrel{\hphantom{=}}
+\pr_1\tilde\sigma_2\otimes\pr_2\tilde\sigma_2|_{\{v_2\}}
+\pr_1\tilde\sigma_2|_{\{v_0,v_1\}}\otimes\pr_2\tilde\sigma_2|_{\{v_1,v_2\}}
+\pr_1\tilde\sigma_2|_{\{v_0\}}\otimes\pr_2\tilde\sigma_2\\
\Sha^*_{(1,1)}(\tau)&=\pr_1\tilde\tau_1\otimes\pr_2\tilde\tau_1|_{\{v_2\}}
+\pr_1\tilde\tau_1|_{\{v_0,v_3\}}\otimes\pr_2\tilde\tau_1|_{\{v_3,v_2\}}
+\pr_1\tilde\tau_1|_{\{v_0\}}\otimes\pr_2\tilde\tau_1\\
&\mathrel{\hphantom{=}}
+\pr_1\tilde\tau_2\otimes\pr_2\tilde\tau_2|_{\{v_2\}}
+\pr_1\tilde\tau_2|_{\{v_0,v_3\}}\otimes\pr_2\tilde\tau_2|_{\{v_3,v_2\}}
+\pr_1\tilde\tau_2|_{\{v_0\}}\otimes\pr_2\tilde\tau_2,
\end{align*}
where \(\pr_i:B_{1,1}\graf\cong B_1\graf\times B_1\graf\to B_1\graf\) is the projection onto the \(i\)-th factor for each \(i=1,2\).

Since \(\pr_i(\tilde\sigma_j),\pr_i(\tilde\tau_j)\in C_2(B_1\graf)\) for each \(i,j\) but \(S_\bullet(\graf)\) has no summand of bigrading \((2,1)\), they are in the kernel of \(g\).
Hence 
\begin{align*}
(g\otimes g)\Sha^*_{(1,1)}(\sigma)&
=g(\pr_1\tilde\sigma_1|_{\{v_0,v_1\}})\otimes g(\pr_2\tilde\sigma_1|_{\{v_1,v_2\}})
+g(\pr_1\tilde\sigma_2|_{\{v_0,v_1\}})\otimes g(\pr_2\tilde\sigma_2|_{\{v_1,v_2\}})\\
(g\otimes g)\Sha^*_{(1,1)}(\tau)&
=g(\pr_1\tilde\tau_1|_{\{v_0,v_1\}})\otimes g(\pr_2\tilde\tau_1|_{\{v_1,v_2\}})
+g(\pr_1\tilde\tau_2|_{\{v_0,v_1\}})\otimes g(\pr_2\tilde\tau_2|_{\{v_1,v_2\}}).
\end{align*}

On the other hand,
\begin{align*}
\pr_1\tilde\sigma_1|_{\{v_0,v_1\}}&=\sfh,&
\pr_2\tilde\sigma_1|_{\{v_1,v_2\}}&=\phi_\sfe(1-t/4)\simeq \sfh',\\
\pr_1\tilde\sigma_2|_{\{v_0,v_1\}}&=\sfh',&
\pr_2\tilde\sigma_2|_{\{v_1,v_2\}}&=\phi_\sfe((2-t)/4),\\
\pr_1\tilde\tau_1|_{\{v_0,v_3\}}&=\sfv,&
\pr_2\tilde\tau_1|_{\{v_3,v_2\}}&=\phi_\sfe(1-(2-s)/4),\\
\pr_1\tilde\tau_2|_{\{v_0,v_3\}}&=\sfh',&
\pr_2\tilde\tau_2|_{\{v_3,v_2\}}&=\phi_\sfe(s/4)\simeq \sfh.
\end{align*}
Moreover, both \(\phi_\sfe((2-t)/4)\) and \(\phi_\sfe(1-(2-s)/4)\) are of bigrading \((1,1)\) and completely contained in the interior of \(\sfe\), and homotopic to the constant map on an interior point of \(\sfe\), whose bigrading is \((0,1)\). 
This mismatch of bigrading implies that their images under \(g\) vanishes as seen in \Cref{example:computation of g} again.

Therefore,
\((g\otimes g)\Sha^*_{(1,1)}(\sigma)=\sfh\otimes\sfh'\) and \((g\otimes g)\Sha^*_{(1,1)}(\tau)=\sfh'\otimes\sfh\), and so
\begin{align*}
\Sha^*_S(f(\sfh\sfh'))&=(\sfh\sfh')\otimes 1 + \sfh\otimes\sfh'-\sfh'\otimes\sfh+1\otimes(\sfh\sfh')
\end{align*}
as desired.

Now for \(\alpha=\sfh_1\dots\sfh_k\), let us assume that \(\sfh\neq\sfh_i\) for all \(i\).
If \(\alpha\) does not involve \(\sfh'\), then \(f(\alpha\sfh)=f(\alpha)\ast f(\sfh)\) by definition.
\begin{align*}
\Sha^*_S(\alpha\sfh)&=(g\otimes g)\Sha^*(f(\alpha\sfh))\\
&=(g\otimes g)\Sha^*(f(\alpha)*f(\sfh))\\
&=(g\otimes g)(\ast\otimes\ast)T(\Sha^*(f(\alpha))\otimes\Sha^*(f(\sfh))).
\end{align*}

Since \(g\) preserves the external product, i.e., \(g\circ *=*\circ (g\otimes g)\), 
\begin{align*}
\Sha^*_S(\alpha\sfh)
&=(\ast\otimes\ast)((g\otimes g)\otimes (g\otimes g)) T(\Sha^*(f(\alpha))\otimes\Sha^*(f(\sfh)))\\
&=(\ast\otimes\ast)T((g\otimes g)\Sha^*(f(\alpha))\otimes (g\otimes g)\Sha^*(f(\sfh)))\\
&=(\ast\otimes\ast)T(\Sha^*_S(f(\alpha))\otimes \Sha^*_S(\sfh))\\
&=(\ast\otimes\ast)T(\Sha^*_S(\alpha)\otimes(\sfh\otimes1+1\otimes\sfh)).
\end{align*}

Then by the induction hypothesis and \Cref{proposition:commutative diagram with external product}, 
\[
\Sha^*_S(\alpha)=\sum_{I\subset[k]}(-1)^{s(I,[k]\setminus I)} \left(\bigwedge_{i\in I}\sfh_i\right)
\otimes\left(\bigwedge_{j\in[k]\setminus I}\sfh_j\right).
\]
If we denote \(\sfh\) by \(\sfh_{k+1}\), we have
\begin{align*}
\Sha^*_S(\alpha\sfh_{k+1})
&=\sum_{I\subset[k]}(-1)^{s(I,[k]\setminus I)}(\ast\otimes\ast) T\left(
\left(\bigwedge_{i\in I}\sfh_i\right)
\otimes \left(\bigwedge_{j\in[k]\setminus I}\sfh_j\right)\otimes
(\sfh_{k+1}\otimes1+1\otimes\sfh_{k+1})
\right)\\
&=\sum_{I\subset[k]}(-1)^{s(I,[k]\setminus I)} (\ast\otimes\ast)\left(
(-1)^{s(\{k+1\},[k]\setminus I)}\left(\bigwedge_{i\in I}\sfh_i\right)\otimes \sfh_{k+1}\otimes 1 \otimes \left(\bigwedge_{j\in [k]\setminus I}\sfh_j\right)\right)\\
&\mathrel{\hphantom{=}}
+\sum_{I\subset[k]}(-1)^{s(I,[k]\setminus I)} (\ast\otimes\ast)\left(1\otimes
\left(\bigwedge_{i\in I}\sfh_i\right)
\otimes \left(\bigwedge_{j\in[k]\setminus I}\sfh_j\right)\otimes \sfh_{k+1}
\right)\\
&=\sum_{I\subset[k]}(-1)^{s(I,[k]\setminus I)+s(\{k+1\},[k]\setminus I)}\left(\bigwedge_{i\in I\cup\{k+1\}}\sfh_i\right)
\otimes \left(\bigwedge_{j\in [k]\setminus I}\sfh_j\right)\\
&\mathrel{\hphantom{=}}
+\sum_{I\subset[k]}(-1)^{s(I,[k]\setminus I)}
\left(\bigwedge_{i\in I}\sfh_i\right)\ast
\left(\bigwedge_{j\in [k+1]\setminus I}\sfh_j\right)\\
&=\sum_{I\subset [k+1]} (-1)^{s(I,[k+1]\setminus I)}\left(\bigwedge_{i\in I}\sfh_i\right)
\otimes \left(\bigwedge_{j\in[k+1]\setminus I}\sfh_j\right).
\end{align*}

Finally, let us denote \(\sfh'\) by \(\sfh_{k+2}\). Then 
\begin{align*}
\Sha^*_S(\alpha\sfh_{k+1}\sfh_{k+2})&=(g\otimes g)\Sha^*(f(\alpha\sfh_{k+1}\sfh_{k+2}))\\
&=(g\otimes g)\Sha^*(f(\alpha)\ast f(\sfh_{k+1}\sfh_{k+2}))\\
&=(g\otimes g)(\ast\otimes\ast) T(\Sha^*(f(\alpha))\otimes\Sha^*(f(\sfh_{k+1}\sfh_{k+2})))\\
&=(\ast\otimes\ast)((g\otimes g)\otimes (g\otimes g)) T(\Sha^*(f(\alpha))\otimes\Sha^*(f(\sfh_{k+1}\sfh_{k+2})))\\
&=(\ast\otimes\ast)T((g\otimes g)\Sha^*(f(\alpha))\otimes (g\otimes g)\Sha^*(f(\sfh_{k+1}\sfh_{k+2})))\\
&=(\ast\otimes\ast)T(\Sha^*_S(\alpha)\otimes\Sha^*_S(\sfh_{k+1}\sfh_{k+2})).
\end{align*}

By the induction hypothesis for \(\alpha, \alpha\sfh_{k+1}\) and the computation of \(\Sha^*_S(\sfh_{k+1}\sfh_{k+2})\) above, we have
\begin{align*}
&\mathrel{\hphantom{=}}\Sha^*_S(\alpha)\otimes \Sha^*_S(\sfh_{k+1}\sfh_{k+2})\\
&=(\ast\otimes\ast)\sum_{I\subset[k]}(-1)^{s(I,[k]\setminus I)+s(\{k+1,k+2\},[k]\setminus I)}
\left(\bigwedge_{i\in I}\sfh_i\right)\otimes \sfh_{k+1}\sfh_{k+2}\otimes
\left(\bigwedge_{j\in [k]\setminus I}\sfh_j\right)\otimes 1
\\
&\mathrel{\hphantom{=}}+
(\ast\otimes\ast)\sum_{I\subset[k]}(-1)^{s(I,[k]\setminus I)+s(\{k+1\},[k]\setminus I)}
\left(\bigwedge_{i\in I}\sfh_i\right)\otimes \sfh_{k+1}\otimes
\left(\bigwedge_{j\in [k]\setminus I}\sfh_j\right)\otimes \sfh_{k+2}
\\
&\mathrel{\hphantom{=}}-
(\ast\otimes\ast)\sum_{I\subset[k]}(-1)^{s(I,[k]\setminus I)+s(\{k+2\},[k]\setminus I)}
\left(\bigwedge_{i\in I}\sfh_i\right)\otimes \sfh_{k+2}\otimes
\left(\bigwedge_{j\in [k]\setminus I}\sfh_j\right)\otimes \sfh_{k+1}
\\
&\mathrel{\hphantom{=}}+
(\ast\otimes\ast)\sum_{I\subset[k]}(-1)^{s(I,[k]\setminus I)}
\left(\bigwedge_{i\in I}\sfh_i\right)\otimes 1\otimes
\left(\bigwedge_{j\in [k]\setminus I}\sfh_j\right)\otimes \sfh_{k+1}\sfh_{k+2}
\end{align*}
and it is the same as follows:
\begin{align*}
\Sha^*_S(\alpha)\otimes \Sha^*_S(\sfh_{k+1}\sfh_{k+2})&=\sum_{I\subset[k]}(-1)^{s(I,[k]\setminus I)+s(\{k+1,k+2\},[k]\setminus I)}
\left(\bigwedge_{i\in I\cup\{k+1,k+2\}}\sfh_i\right)\otimes
\left(\bigwedge_{j\in [k]\setminus I}\sfh_j\right)
\\
&\mathrel{\hphantom{=}}+
\sum_{I\subset[k]}(-1)^{s(I,[k]\setminus I)+s(\{k+1\},[k]\setminus I)}
\left(\bigwedge_{i\in I\cup\{k+1\}}\sfh_i\right)\otimes
\left(\bigwedge_{j\in [k]\cup\{k+2\}\setminus I}\sfh_j\right)
\\
&\mathrel{\hphantom{=}}-
\sum_{I\subset[k]}(-1)^{s(I,[k]\setminus I)+s(\{k+2\},[k]\setminus I)}
\left(\bigwedge_{i\in I\cup\{k+2\}}\sfh_i\right)\otimes
\left(\bigwedge_{j\in [k]\cup\{k+1\}\setminus I}\sfh_j\right)
\\
&\mathrel{\hphantom{=}}+
\sum_{I\subset[k]}(-1)^{s(I,[k]\setminus I)}
\left(\bigwedge_{i\in I}\sfh_i\right)\otimes
\left(\bigwedge_{j\in [k]\cup\{k+1,k+2\}\setminus I}\sfh_j\right)\\
&=\sum_{I\subset[k+2]}(-1)^{s(I,[k+2]\setminus I)}
\left(\bigwedge_{i\in I}\sfh_i\right)\otimes
\left(\bigwedge_{j\in [k+2]}\sfh_j\right).\qedhere
\end{align*}
\end{proof}

\subsection{Comultiplication \(\sha^*\) on \(\ring[E]\)}
Let \(\bbZ_{\ge0}^E\) be the set of sequences \(\bfa=(a_\sfe)_{\sfe\in E}\) of nonnegative integers indexed by the set \(E\), and let \(\|\bfa\|=\sum_{\sfe\in E}a_\sfe\).
For \(\bfa=(a_\sfe)_{\sfe\in E},\bfb=(b_\sfe)_{\sfe\in E}\in\bbZ_{\ge0}^E\), we write \(\bfb\le\bfa\) if \(b_\sfe\le a_\sfe\) for any \(\sfe\in E\).
For each \(\bfb\in \bbZ_{\ge0}^e\), we also define \(\binom{\bfa}{\bfb}=\prod_{\sfe\in E}\binom{a_\sfe}{b_\sfe}\). Then \(\binom{\bfa}{\bfb}\) is nonzero if and only if \(\bfb\le\bfa\).
For a mononomial \(m\in \ring[E]\), we say that \(m\) is \emph{\(\sfe\)-free} if \(m\) has no \(\sfe\), or equivalently, \(a_\bfe=0\) if \(m=E^\bfa\).

Let \(E'=\{\sfe':\sfe\in E\}\) be a copy of \(E\).
Then there is a canonical isomorphism \(F:\ring[E,E']\to\ring[E]\otimes\ring[E]\), where \(\sfe\) and \(\sfe'\) correspond to \(\sfe\otimes 1\) and \(1\otimes\sfe\), respectively.
Since \(\sha^*:\ring[E]\to\ring[E]\otimes\ring[E]\) is an algebra homomorphism, we have the following commutative diagram:
\[
\begin{tikzcd}[column sep=3pc]
\ring[E] \ar[rr,"\sha^*"] \ar[rd,"\Phi"'] && \ring[E]\otimes\ring[E]\\
& \ring[E,E'] \ar[ru, "\cong","F"']
\end{tikzcd}
\]
where \(\Phi\) is defined as for each \(p(E)\in\ring[E]\),
\[
\Phi(p(E))= p(\sfe+\sfe':\sfe\in E).
\]

\begin{example}
For \(p(\sfe_1,\sfe_2)=\sfe_1\sfe_2+\sfe_1\), we have
\begin{align*}
\sha^*(p(\sfe_1,\sfe_2)) &= (\sfe_1\otimes1+1\otimes\sfe_2)(\sfe_2\otimes1+1\otimes\sfe_2)+(\sfe_1\otimes1+1\otimes\sfe_1)\\
\Phi(p(\sfe_1,\sfe_2)) &= p(\sfe_1+\sfe_1', \sfe_2+\sfe_2')\\
&=(\sfe_1+\sfe_1')(\sfe_2+\sfe_2')+(\sfe_1+\sfe_1').
\end{align*}
Then under the isomorphism \(\ring[\sfe_1,\sfe_2,\sfe_1',\sfe_2']\cong\ring[\sfe_1,\sfe_2]\otimes\ring[\sfe_1,\sfe_2]\), the polynomial \(\Phi(p(\sfe_1,\sfe_2))\) will be mapped to \(\sha^*(p(\sfe_1,\sfe_2))\).
\end{example}

\begin{lemma}\label{lemma:comultiplication on ring}
Let \(\bfa\in \bbZ^E_{\ge0}\) be a sequence of nonnegative integers. Then
\begin{align*}
\sha^*(E^\sfa) &=\sum_{\bfb\le\bfa}\binom{\bfa}{\bfb}E^\bfb\otimes E^{\bfa-\bfb}.
\end{align*}
\end{lemma}
\begin{proof}
As seen above, \(\sha^*(E^\sfa)=F(\Phi(E^\sfa))\). Hence 
\begin{align*}
\sha^*(E^\sfa)&=F\left(\prod_{\sfe\in E}(\sfe+\sfe')^{a_\sfe}\right)
=\sum_{\bfb\le\bfa}\binom{\bfa}{\bfb}F(\sfe^{b_\sfe}\sfe'^{a_\sfe-b_\sfe})
=\sum_{\bfb\le\bfa}\binom{\bfa}{\bfb}E^\sfb\otimes E^{\sfa-\sfb}.\qedhere
\end{align*}
\end{proof}

\begin{definition}[{The subalgebra \(\ring_0[E]\)}]\label{definition:R0}
We define a subalgebra \(\ring_0[E]\) of \(\ring[E]\) generated by elements of the form \(\sfe_i-\sfe_j\) for \(\sfe_i,\sfe_j\in E\), or equivalently, of the form \(\sfe-\sfe_0\) for a chosen edge \(\sfe_0\) and \(\sfe\in E\).
\end{definition}

For the latter purpose, we consider the following situation.
For a chosen edge \(\sfe_0\in E\), let \(\pi:\ring[E]\otimes\ring[E]\to\ring[E]\otimes\ring[\sfe_0]\) be an algebra homomorphism that maps \(\sfe\otimes1\) and \(1\otimes\sfe\) to \(\sfe\otimes1\) and \(1\otimes\sfe_0\) for every \(\sfe\in E\), respectively.

\begin{proposition}\label{proposition:graded augmentation}
Under the composition \(\ring[E]\stackrel{\sha^*}{\to}\ring[E]\otimes\ring[E]\stackrel{\pi}{\to}\ring[E]\otimes\ring[\sfe_0]\),
a polynomial \(p(E)\) has the image \(p(E)\otimes 1\) if and only if \(p(E)\in \ring_0[E]\).
\end{proposition}
\begin{proof}
Let us consider an algebra isomorphism \(F:\ring[E,\sfe_0']\to\ring[E]\otimes\ring[\sfe_0]\) that maps \(\sfe\) and \(\sfe_0'\) to \(\sfe\otimes 1\) and \(1\otimes\sfe\) for every \(\sfe\in E\).
Then for each polynomial \(p(E)\), it is easy to observe that
\[
\pi\circ\sha^*(p(E)) = F(p(\sfe+\sfe_0':\sfe\in E)).
\]

Let \(\bar\sfe = \sfe-\sfe_0\) for each \(\sfe_i\in E\) and \(\bar E=\{\bar\sfe:\sfe_0\neq\sfe\in E\}\).
Then there is another polynomial \(q(\sfe_0,\bar E)\) of \(\sfe_0\) and elements in \(\bar E\) such that \(p(E)=q(\sfe_0,\bar E)\) in \(\ring[E]\).
Notice that under the evaluation \(\sfe\mapsto\sfe+\sfe_0'\) for every \(\sfe\in E\), each variable \(\bar\sfe=\sfe-\sfe_0=(\sfe+\sfe_0')-(\sfe_0+\sfe_0')\) remains unchanged and so
\[
\pi\circ\sha^*(p(E)) = f(p(\sfe+\sfe_0':\sfe\in E)) = f(q(\sfe_0+\sfe_0',\bar E)).
\]

Therefore \(\pi\circ\sha^*(p(E))=p(E)\otimes 1\) if and only if \(q(\sfe_0+\sfe_0',\bar E)\) is constant with respect to \(\sfe_0'\), which is also equivalent to that \(q(\sfe_0,\bar E)\) is constant with respect to \(\sfe_0\).
In other words, \(p(E) = q(0,\bar E)\) is a polynomial of variables in \(\bar E\), which is equivalent to \(p(E)=q(0,\bar E)\in \ring_0[E]\) as desired.
\end{proof}

\section{Applications to bunches of grapes}
\subsection{Bunches of grapes}
Let \(\gamma=(\sfv_0,\dots,\sfv_n)\) be a simple cycle of a connected graph \(\graf\),
that is, \(\sfv_i\) and \(\sfv_{i+1}\) are joined by an edge for each
\(0\le i<n\), \(\sfv_0=\sfv_n\), and \(\sfv_1,\dots, \sfv_n\) are
distinct.
The number \(n\) is called the \emph{length} of \(\gamma\).

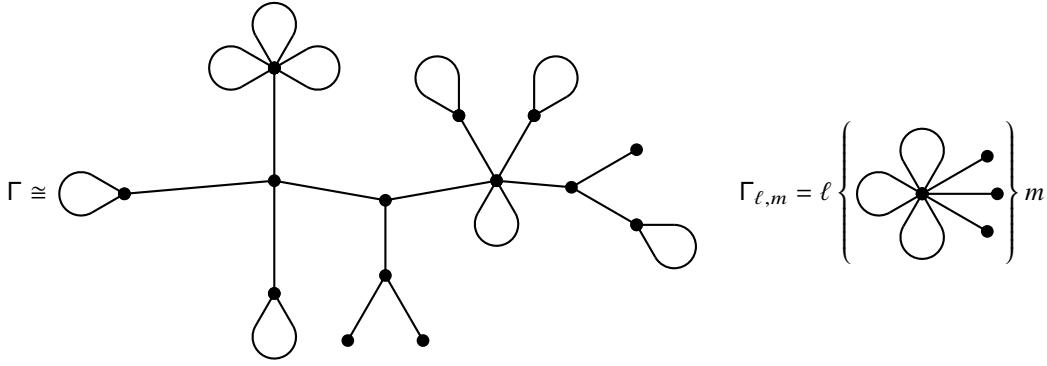
\begin{figure}[ht]
\begin{align*}
\graf&\cong\begin{tikzpicture}[baseline=-.5ex]
\draw[thick,fill] (0,0) circle (2pt) node (A) {} -- ++(5:2) circle(2pt) node(B) {} -- +(0,1.5) circle (2pt) node (C) {} +(0,0) -- +(0,-1.5) circle (2pt) node(D) {} +(0,0) -- ++(-10:1.5) circle(2pt) node(E) {} -- ++(10:1.5) circle(2pt) node(F) {} -- +(120:1) circle(2pt) node(G) {} +(0,0) -- +(60:1) circle (2pt) node(H) {} +(0,0) -- ++(-5:1) circle (2pt) node(I) {} -- +(30:1) circle(2pt) node(J) {} +(0,0) -- +(-30:1) circle (2pt) node(K) {};
\draw[thick, fill] (E.center) -- ++(0,-1) circle (2pt) -- +(-60:1) circle (2pt) +(0,0) -- +(-120:1) circle (2pt);
\grape[180]{A};
\grape[0]{C}; \grape[90]{C}; \grape[180]{C};
\grape[-90]{D};
\grape[-90]{F};
\grape[120]{G};
\grape[60]{H};
\grape[-30]{K};
\end{tikzpicture}&
\graf_{\ell,m}&=\ell\left\{\begin{tikzpicture}[baseline=-.5ex]
\draw[thick,fill] (0,0) circle (2pt) node (X) {} -- (0:1) circle (2pt) (0,0) -- (30:1) circle (2pt) (0,0) -- (-30:1) circle (2pt);
\grape[90]{X};
\grape[180]{X};
\grape[270]{X};
\end{tikzpicture}
\right\}m
\end{align*}
\caption{A bunch of grapes (left) and an elementary bunch of grapes (right).}
\label{figure:example of a bunch of grapes}
\end{figure}

The \emph{circumference} \(\sfC(\graf)\) of \(\graf\) is the maximum
number of edges in a simple cycle of \(\graf\). The \emph{topological
 circumference} \(\TC(\graf)=\min\{\sfC(\graf'): \graf\cong\graf'\}\) of \(\graf\) is the minimum
circumference among all graphs homeomorphic to \(\graf\).
Then \(\TC(\graf)=0\) if and only if \(\graf\) is a tree, and any graph \(\graf\) with \(\TC(\graf)=1\) is obtained by attaching loops to a
tree. 
We call a graph with topological circumference \(1\) a \emph{bunch of grapes}.
In particular, a bunch of grapes with at most one essential vertex will be said to be \emph{elementary}, which can be completely classified up to homeomorphism by the numbers of leaves \(m\) and loops \(\ell\) and denoted by \(\graf_{\ell,m}\).
See \Cref{figure:example of a bunch of grapes}.

\begin{theorem}\cite{KP2012}\label{theorem:basis for elementary}
For each \(k\ge 0\), the \(i\)-th homology group \(H_i(B_k(\graf_{\ell,m}))\) is a free \(\ring\)-module of rank \(1\) if \(i=0\),
\begin{equation}\label{eq:first homology}
N_{\ell,m}(k)\coloneqq (2\ell+m-2)\binom {k+\ell+m-2}{\ell+m-1} - \binom{k+\ell+m-2}{\ell+m-2}+1.
\end{equation}
if \(i=1\),
or \(0\) if \(i\ge 2\).
\end{theorem}

In general, any bunch of grapes \(\graf\) can be uniquely determined up to homeomorphism by a pair of tree \(\sfT\) and a function \(\loops:V(\sfT)\to\bbZ_{\ge0}\), where \(\loops(\sfv)\) is the number of cycles attached to a given vertex \(\sfv\in V(\sfT)\).
In this case, we write \(\graf\cong(\sfT,\loops)\) and call \(\sfT\) the \emph{stem} of \(\graf\).
For example, the stem and the function \(\loops\) for a bunch of grapes in \Cref{figure:example of a bunch of grapes} is as follows:
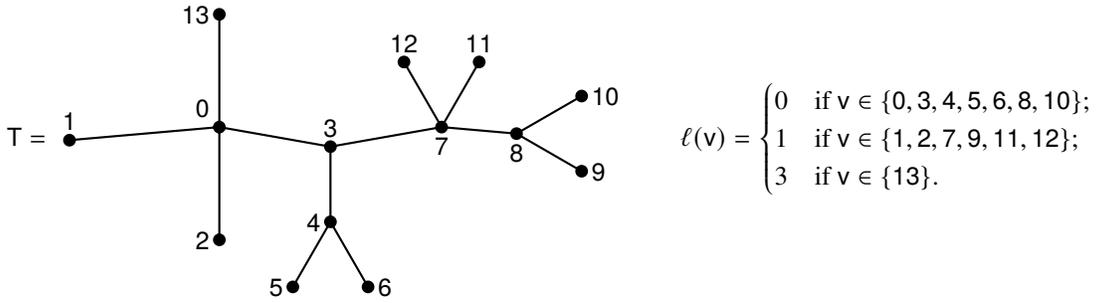
\[
\sfT=\begin{tikzpicture}[baseline=-.5ex]
\draw[thick,fill] (0,0) circle (2pt) node (A) {} node[above] {\(\mathsf{1}\)} -- ++(5:2) circle(2pt) node(B) {} node[above left] {\(\mathsf{0}\)} -- +(0,1.5) circle (2pt) node (C) {} node[left] {\(\mathsf{13}\)} +(0,0) -- +(0,-1.5) circle (2pt) node(D) {} node[left] {\(\mathsf{2}\)} +(0,0) -- ++(-10:1.5) circle(2pt) node(E) {} node[above] {\(\mathsf{3}\)} -- ++(10:1.5) circle(2pt) node(F) {} node[below] {\(\mathsf{7}\)} -- +(120:1) circle(2pt) node(G) {} node[above] {\(\mathsf{12}\)} +(0,0) -- +(60:1) circle (2pt) node(H) {} node[above] {\(\mathsf{11}\)} +(0,0) -- ++(-5:1) circle (2pt) node(I) {} node[below] {\(\mathsf{8}\)} -- +(30:1) circle(2pt) node(J) {} node[right] {\(\mathsf{10}\)} +(0,0) -- +(-30:1) circle (2pt) node(K) {} node[right] {\(\mathsf{9}\)};
\draw[thick, fill] (E.center) -- ++(0,-1) circle (2pt) node[left] {\(\mathsf{4}\)} -- +(-60:1) circle (2pt) node[right] {\(\mathsf{6}\)} +(0,0) -- +(-120:1) circle (2pt) node[left] {\(\mathsf{5}\)};
\end{tikzpicture}\qquad
\loops(\sfv)=\begin{cases}
0 & \mbox{if }\sfv\in \{\mathsf{0}, \mathsf{3}, \mathsf{4}, \mathsf{5}, \mathsf{6}, \mathsf{8}, \mathsf{10}\};\\
1 & \mbox{if }\sfv\in \{\mathsf{1}, \mathsf{2}, \mathsf{7}, \mathsf{9}, \mathsf{11}, \mathsf{12}\};\\
3 & \mbox{if }\sfv\in\{\mathsf{13}\}.
\end{cases}
\]
\begin{assumption}
From now on, we assume that every bunch of grapes \(\graf=(\sfT,\loops)\) has at least one essential vertex and no bivalent vertices.
\end{assumption}

We consider the general case first, when the stem \(\sfT\) has at least one edge.
Then we fix a planar embedding of \(\graf\) to \(\bbR^2\) without nested cycles, namely, the interiors of discs that each cycle bounds are disjoint.
We also fix a pair \((\sfv_0,\sfe_0)\) of
a vertex \( \sfv_0 \) and an edge \(\sfe_0\) of \( \sfT \)
adjacent to \(\sfv_0\), the \emph{root vertex} and the \emph{root edge}, respectively.
We call this pair an \emph{oriented root} of \(\graf\).

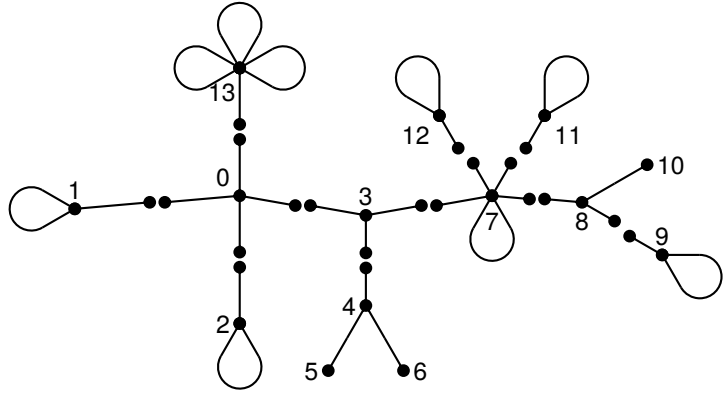
\begin{figure}[ht]
\[
\begin{tikzpicture}[baseline=-.5ex]
\draw[thick,fill] (0,0) circle (2pt) node (A) {} node[above] {\(\mathsf{1}\)} -- ++(5:1) circle (2pt);
\grape[180]{A};
\begin{scope}[xshift=0.2cm]
	\draw[thick,fill] (0,0) +(5:1) circle (2pt) -- ++(5:2) circle(2pt) node(B) {} node[above left] {\(\mathsf{0}\)} -- +(0,0.75) circle (2pt) +(0,0) -- +(0,-0.75) circle (2pt) +(0,0) -- ++(-10:0.75) circle (2pt);

	\begin{scope}[yshift=0.2cm]
		\draw[thick,fill] (0,0) ++(5:2) +(0,0.75) circle (2pt) -- +(0,1.5) circle (2pt) node (C) {} node[below=2ex, left=-.5ex] {\(\mathsf{13}\)};
		\grape[0]{C}; \grape[90]{C}; \grape[180]{C};
	\end{scope}

	\begin{scope}[yshift=-0.2cm]
		\draw[thick,fill] (0,0) ++(5:2) +(0,-0.75) circle (2pt) -- +(0,-1.5) circle (2pt) node (D) {} node[left] {\(\mathsf{2}\)};
		\grape[-90]{D};
	\end{scope}

	\begin{scope}[xshift=0.2cm]
		\draw[thick,fill] (0,0) ++(5:2) ++(-10:0.75) circle (2pt) -- ++(-10:0.75) circle(2pt) node(E) {} node[above] {\(\mathsf{3}\)} -- +(10:0.75) circle (2pt) +(0,0) -- +(0,-0.5) circle (2pt);

		\begin{scope}[yshift=-0.2cm]
			\draw[thick,fill] (0,0) ++(5:2) ++(-10:1.5) ++(0,-0.5) circle (2pt) -- ++(0,-0.5) circle (2pt) node[left] {\(\mathsf{4}\)} -- +(-60:0.5) +(0,0) -- +(-120:0.5);

			\begin{scope}[yshift=0cm]
				\draw[thick,fill] (0,0) ++(5:2) ++(-10:1.5) ++(0,-1) +(-120:0.5) -- +(-120:1) circle (2pt) node[left] {\(\mathsf{5}\)} +(-60:0.5) -- +(-60:1) circle (2pt) node[right] {\(\mathsf{6}\)};
			\end{scope}
		\end{scope}

		\begin{scope}[xshift=0.2cm]
			\draw[thick,fill] (0,0) ++(5:2) ++(-10:1.5) ++(10:0.75) circle(2pt) -- ++(10:0.75) circle(2pt) node(F) {} node[below=0.5ex] {\(\mathsf{7}\)} -- +(120:0.5) circle (2pt) +(0,0) -- +(60:0.5) circle (2pt) +(0,0) -- +(-5:0.5) circle (2pt);
			\grape[-90]{F};

			\begin{scope}[xshift=-0.2cm, yshift=0.2cm]
				\draw[thick,fill] (0,0) ++(5:2) ++(-10:1.5) ++(10:1.5) +(120:0.5) circle (2pt) -- +(120:1) circle (2pt) node(G) {} node[below left] {\(\mathsf{12}\)};
				\grape[120]{G};
			\end{scope}
			\begin{scope}[xshift=0.2cm, yshift=0.2cm]
				\draw[thick,fill] (0,0) ++(5:2) ++(-10:1.5) ++(10:1.5) +(60:0.5) circle (2pt) -- +(60:1) circle (2pt) node(H) {} node[below right] {\(\mathsf{11}\)};
				\grape[60]{H};
			\end{scope}

			\begin{scope}[xshift=0.2cm]
				\draw[thick,fill] (0,0) ++(5:2) ++(-10:1.5) ++(10:1.5) ++(-5:0.5) circle (2pt) -- ++(-5:0.5) circle (2pt) node(I) {} node[below] {\(\mathsf{8}\)} -- +(30:0.5) +(0,0) -- +(-30:0.5) circle (2pt);

				\begin{scope}[xshift=0cm]
					\draw[thick,fill] (0,0) ++(5:2) ++(-10:1.5) ++(10:1.5) ++(-5:1) +(30:0.5) -- +(30:1) circle (2pt) node[right] {\(\mathsf{10}\)};
				\end{scope}
				\begin{scope}[xshift=0.2cm,yshift=-.2cm]
					\draw[thick,fill] (0,0) ++(5:2) ++(-10:1.5) ++(10:1.5) ++(-5:1) +(-30:0.5) circle (2pt) -- +(-30:1) circle (2pt) node(K) {} node[above] {\(\mathsf{9}\)};
					\grape[-30]{K};
				\end{scope}
			\end{scope}
		\end{scope}
	\end{scope}
\end{scope}

\end{tikzpicture}
\]
\caption{The decomposition of \(\graf\) along the stem.}
\label{figure:decomposion along the stem}
\end{figure}

For each vertex \(\sfv\neq\sfv_0\) in \(\sfT\), there is a unique edge \(\sfe^\sfv\) adjacent to \(\sfv\), which is the first edge of the shortest path from \(\sfv\) to \(\sfv_0\).
Then for each essential vertex \(\sfv\), by cutting the edge \(\sfe^\sfv\) into two pieces, we obtain a collection of elementary bunches of grapes \(\{\graf^\sfv=\graf_{\ell(\sfv),m(\sfv)}:\sfv\in V^\ess(\graf)\}\), where \(m(\sfv)=\val_\sfT(\sfv)\). 
We call each \(\graf_{\ell(\sfv),m(\sfv)}\) the \emph{local graph} of \(\graf\) at \(\sfv\).
See \Cref{figure:decomposion along the stem}. The oriented edge is denoted by a red arrow.

Now we label edges in the local graph \(\graf^\sfv\) as follows:
At first, there is a unique edge in \(\graf^\sfv\), which is a part of the edge \(\sfe^\sfv\) in \(\graf\) and will be denoted by \(\sfe^\sfv_1\).
By using the planar embedding of \(\graf\) that induces a planar embedding of \(\graf^\sfv\), we label all other edges with \(\sfe^\sfv_2,\dots,\sfe^\sfv_{2\ell(\sfv)+m(\sfv)}\) in a counterclockwise manner starting from \(\sfe^\sfv_1\).
We also label half-edges adjacent to \(\sfv\) with \(\sfh^\sfv_1,\dots,\sfh^\sfv_{2\ell(\sfv)+m(\sfv)}\) in the same manner so that \(\sfh^\sfv_1\subset\sfe^\sfv_1\).
For \(i\ge 2\), we call a half-edge \(\sfh^\sfv_i\) contained in \(\sfe\) the \emph{second half-edge} of \(\sfe\) if \(\sfh^\sfv_{i-1}\cup\sfh^\sfv_i=\sfe\), or the \emph{first half-edge} of \(\sfe\) otherwise.

Finally, we label each loop adjacent to \(\sfv\) in \(\graf\) as the label of corresponding cycle in a local graph \(\graf^\sfv\). See \Cref{figure:labelling1,figure:labelling2}. We omit some labels for brevity.

\begin{figure}[ht]
  \subcaptionbox{Labels on edges of \(\Gamma\).
    \label{figure:labelling1}}[.8\textwidth]{
\(
\begin{tikzpicture}[baseline=-.5ex]
\draw[thick,fill] (0,0) circle (2pt) node (A) {} node[above] {\(\mathsf{1}\)} -- ++(5:2) node[midway, above] {\(\sfe^{\mathsf{1}}\)} circle(2pt) node(B) {} node[above left] {\(\mathsf{0}\)};
\draw[thick,fill] (B.center) +(0,1.5) circle (2pt) node (C) {} node[left=0.5ex] {\(\mathsf{13}\)};
\draw[very thick,fill,<-, red] (B.center)+(0,2pt) -- (C.center) node[midway, right] {\(\sfe^{\mathsf{13}}\)};
\draw[thick,fill] (B.center) -- +(0,-1.5) node[midway, right] {\(\sfe^{\mathsf{2}}\)} circle (2pt) node(D) {} node[left] {\(\mathsf{2}\)} +(0,0) -- ++(-10:1.5) node[midway, above] {\(\sfe^{\mathsf{3}}\)} circle(2pt) node(E) {} node[above] {\(\mathsf{3}\)} -- ++(10:1.5) node[midway, below] {\(\sfe^{\mathsf{7}}\)} circle(2pt) node(F) {} node[below=0.5ex] {\(\mathsf{7}\)} -- +(120:1) node[midway, left] {\(\sfe^{\mathsf{12}}\)} circle(2pt) node(G) {} node[xshift=-1.2ex,yshift=2ex] {\(\mathsf{12}\)} +(0,0) -- +(60:1) node[midway, right] {\(\sfe^{\mathsf{11}}\)} circle (2pt) node(H) {} node[xshift=1.2ex,yshift=2ex] {\(\mathsf{11}\)} +(0,0) -- ++(-5:1) node[midway, below] {\(\sfe^{\mathsf{8}}\)} circle (2pt) node(I) {} node[below] {\(\mathsf{8}\)} -- +(30:1) node[midway, above] {\(\sfe^{\mathsf{10}}\)} circle(2pt) node(J) {} node[right] {\(\mathsf{10}\)} +(0,0) -- +(-30:1) node[midway, below] {\(\sfe^{\mathsf{9}}\)} circle (2pt) node(K) {} node[above] {\(\mathsf{9}\)};
\draw[thick, fill] (E.center) -- ++(0,-1) node[midway, right] {\(\sfe^{\mathsf{4}}\)} circle (2pt) node[left] {\(\mathsf{4}\)} -- +(-60:1) node[midway, right] {\(\sfe^{\mathsf{6}}\)} circle (2pt) node[right] {\(\mathsf{6}\)} +(0,0) -- +(-120:1) node[midway, left] {\(\sfe^{\mathsf{5}}\)} circle (2pt) node[left] {\(\mathsf{5}\)};
\grape[180]{A};
\draw (A.center)+(180:{1.5/sqrt(3)}) node[left] {\(\sfe_2^{\mathsf{1}}\)};
\grape[0]{C}; \grape[90]{C}; \grape[180]{C};
\draw (C.center)+(0:{1.5/sqrt(3)}) node[right] {\(\sfe_2^{\mathsf{13}}\)};
\draw (C.center)+(90:{1.5/sqrt(3)}) node[above] {\(\sfe_3^{\mathsf{13}}\)};
\draw (C.center)+(180:{1.5/sqrt(3)}) node[left] {\(\sfe_4^{\mathsf{13}}\)};
\grape[-90]{D};
\draw (D.center)+(-90:{1.5/sqrt(3)}) node[below] {\(\sfe_2^{\mathsf{2}}\)};
\grape[-90]{F};
\draw (F.center)+(-90:{1.5/sqrt(3)}) node[below] {\(\sfe_2^{\mathsf{7}}\)};
\grape[120]{G};
\draw (G.center)+(120:{1.5/sqrt(3)}) node[above] {\(\sfe_2^{\mathsf{12}}\)};
\grape[60]{H};
\draw (H.center)+(60:{1.5/sqrt(3)}) node[above] {\(\sfe_2^{\mathsf{11}}\)};
\grape[-30]{K};
\draw (K.center)+(-30:{1.5/sqrt(3)}) node[right] {\(\sfe_2^{\mathsf{9}}\)};
\end{tikzpicture}
\)
}
\subcaptionbox{Labels on edges of the decomposition of \(\Gamma\).
  \label{figure:labelling2}}[.8\textwidth]{
\(
\begin{tikzpicture}[baseline=-.5ex]
\draw[thick,fill] (0,0) circle (2pt) node (A) {} node[above] {\(\mathsf{1}\)} [blue] -- ++(5:1) circle (2pt) node[midway,above] {\(\sfe_1^{\mathsf{1}}\)};
\draw (-180:{1.5/sqrt(3)}) node[left] {\(\sfe_2^{\mathsf{1}}\)};
\grape[180]{A};
\begin{scope}[xshift=0.2cm]
	\draw[thick,fill] (0,0) +(5:1) circle (2pt) -- ++(5:2) circle(2pt) node(B) {} node[above left] {\(\mathsf{0}\)} +(0,0) -- +(0,-0.75) circle (2pt) +(0,0) -- +(-10:0.75) circle (2pt);
	\draw[very thick, fill, red,<-] (B.center)+(0,2pt) -- +(0,0.75) node[midway, right] {\(\sfe_1^{\mathsf{0}}\)};
	\draw[thick, fill, blue] (B.center) +(0,0.75) circle (2pt);

	\begin{scope}[yshift=0.2cm]
		\draw[thick,fill, blue] (0,0) ++(5:2) +(0,0.75) circle (2pt) -- +(0,1.5) node[near start,right] {\(\sfe_1^{\mathsf{13}}\)};
		\draw[thick,fill] (0,0) ++(5:2) +(0,0.75) +(0,1.5) circle (2pt) node (C) {}  node[left=0.5ex] {\(\mathsf{13}\)};
		\grape[0]{C}; \grape[90]{C}; \grape[180]{C};
		\draw (C.center) +(0:{1.5/sqrt(3)}) node[right] {\(\sfe_2^{\mathsf{13}}\)};
		\draw (C.center) +(90:{1.5/sqrt(3)}) node[above] {\(\sfe_3^{\mathsf{13}}\)};
		\draw (C.center) +(180:{1.5/sqrt(3)}) node[left] {\(\sfe_4^{\mathsf{13}}\)};
	\end{scope}

	\begin{scope}[yshift=-0.2cm]
		\draw[thick,fill, blue] (0,0) ++(5:2) +(0,-0.75) circle (2pt) -- +(0,-1.5) node[midway,right] {\(\sfe_1^{\mathsf{2}}\)};
		\draw[thick,fill] (0,0) ++(5:2) +(0,-0.75) +(0,-1.5) circle (2pt) node (D) {} node[left] {\(\mathsf{2}\)};
		\grape[-90]{D};
		\draw (D.center) +(-90:{1.5/sqrt(3)}) node[below] {\(\sfe_2^{\mathsf{2}}\)};
	\end{scope}

	\begin{scope}[xshift=0.2cm]
		\draw[thick,fill, blue] (0,0) ++(5:2) ++(-10:0.75) circle (2pt) -- ++(-10:0.75) node[midway,above] {\(\sfe_1^{\mathsf{3}}\)};
		\draw[thick,fill] (0,0) ++(5:2) ++(-10:0.75) ++(-10:0.75) circle(2pt) node(E) {} node[above] {\(\mathsf{3}\)} -- +(10:0.75) circle (2pt) +(0,0) -- +(0,-0.5) circle (2pt);

		\begin{scope}[yshift=-0.2cm]
			\draw[thick,fill, blue] (0,0) ++(5:2) ++(-10:1.5) ++(0,-0.5) circle (2pt) -- ++(0,-0.5) node[midway,right] {\(\sfe_1^{\mathsf{4}}\)};
			\draw[thick,fill] (0,0) ++(5:2) ++(-10:1.5) ++(0,-0.5) ++(0,-0.5) circle (2pt) node[left] {\(\mathsf{4}\)} -- +(-60:0.5) +(0,0) -- +(-120:0.5);

			\begin{scope}[yshift=0cm]
				\draw[thick,fill] (0,0) ++(5:2) ++(-10:1.5) ++(0,-1) +(-120:0.5) -- +(-120:1) circle (2pt) node[left] {\(\mathsf{5}\)} +(-60:0.5) -- +(-60:1) circle (2pt) node[right] {\(\mathsf{6}\)};
			\end{scope}
		\end{scope}

		\begin{scope}[xshift=0.2cm]
			\draw[thick,fill,blue] (0,0) ++(5:2) ++(-10:1.5) ++(10:0.75) circle(2pt) -- ++(10:0.75) node[near start,below] {\(\sfe_1^{\mathsf{7}}\)};
			\draw[thick,fill] (0,0) ++(5:2) ++(-10:1.5) ++(10:0.75) ++(10:0.75) circle(2pt) node(F) {} node[below=0.5ex] {\(\mathsf{7}\)} -- +(120:0.5) circle (2pt) +(0,0) -- +(60:0.5) circle (2pt) +(0,0) -- +(-5:0.5) circle (2pt);
			\grape[-90]{F};
			\draw (F.center) +(-90:{1.5/sqrt(3)}) node[below] {\(\sfe_2^{\mathsf{7}}\)};

			\begin{scope}[xshift=-.2cm, yshift=0.2cm]
				\draw[thick,fill,blue] (0,0) ++(5:2) ++(-10:1.5) ++(10:1.5) +(120:0.5) circle (2pt) -- +(120:1) node[near start,left] {\(\sfe_1^{\mathsf{12}}\)};
				\draw[thick,fill] (0,0) ++(5:2) ++(-10:1.5) ++(10:1.5) +(120:1) circle (2pt) node(G) {} node[xshift=-1.2ex,yshift=2ex] {\(\mathsf{12}\)};
				\grape[120]{G};
				\draw (G.center) +(120:{1.5/sqrt(3)}) node[above] {\(\sfe_2^{\mathsf{12}}\)};
			\end{scope}
			\begin{scope}[xshift=.2cm, yshift=0.2cm]
				\draw[thick,fill,blue] (0,0) ++(5:2) ++(-10:1.5) ++(10:1.5) +(60:0.5) circle (2pt) -- +(60:1) node[near start,right] {\(\sfe_1^{\mathsf{11}}\)};
				\draw[thick,fill] (0,0) ++(5:2) ++(-10:1.5) ++(10:1.5) +(60:1) node(H) {} node[xshift=1.2ex,yshift=2ex] {\(\mathsf{11}\)};
				\grape[60]{H};
				\draw (H.center) +(60:{1.5/sqrt(3)}) node[above] {\(\sfe_2^{\mathsf{11}}\)};
			\end{scope}

			\begin{scope}[xshift=0.2cm]
				\draw[thick,fill,blue] (0,0) ++(5:2) ++(-10:1.5) ++(10:1.5) ++(-5:0.5) circle (2pt) -- ++(-5:0.5) node[midway, above] {\(\sfe_1^{\mathsf{8}}\)};
				\draw[thick,fill] (0,0) ++(5:2) ++(-10:1.5) ++(10:1.5) ++(-5:0.5) ++(-5:0.5) circle (2pt) node(I) {} node[below] {\(\mathsf{8}\)} -- +(30:0.5) +(0,0) -- +(-30:0.5) circle (2pt);

				\begin{scope}[xshift=0cm]
					\draw[thick,fill] (0,0) ++(5:2) ++(-10:1.5) ++(10:1.5) ++(-5:1) +(30:0.5) -- +(30:1) circle (2pt) node[right] {\(\mathsf{10}\)};
				\end{scope}
				\begin{scope}[xshift=0.2cm,yshift=-.2cm]
					\draw[thick,fill,blue] (0,0) ++(5:2) ++(-10:1.5) ++(10:1.5) ++(-5:1) +(-30:0.5) circle (2pt) -- +(-30:1) node[midway,below left] {\(\sfe_1^{\mathsf{9}}\)};
					\draw[thick,fill] (0,0) ++(5:2) ++(-10:1.5) ++(10:1.5) ++(-5:1) +(-30:0.5) +(-30:1) circle (2pt) node(K) {} node[above] {\(\mathsf{9}\)};
					\grape[-30]{K};
					\draw (K.center) +(-30:{1.5/sqrt(3)}) node[right] {\(\sfe_2^{\mathsf{9}}\)};
				\end{scope}
			\end{scope}
		\end{scope}
	\end{scope}
\end{scope}
\end{tikzpicture}
\)
}
\caption{Labeling convention of edges in \(\Gamma\) and its decomposition.}
\label{figure:labelling}
\end{figure}

For a root edge \(\sfe_0\), let \(\iota_0:\sfe_0\to\sfe_0\) be an embedding, whose image is disjoint from the endpoints of \(\sfe_0\).
For each \(\sfv\in V^\ess(\graf)\), by shrinking edges of local graphs a little if necessary, we can consider embeddings \(\iota^\sfv:\graf^\sfv \to \graf\) and \(\iota_0:\sfe_0\to \sfe_0\) such that the images of \(\iota^\sfv\) are pairwise disjoint for all \(\sfv\in V^\ess(\graf)\), and disjoint from the image of \(\iota_0\).
Therefore these embeddings can be used to define the external product 
\[
\ast: H_\bullet(B\sfe_0)\otimes\bigotimes_{\sfv\in V^\ess(\graf)} H_\bullet(B\graf^\sfv) \to H_\bullet(B\graf).
\]
It is easy to check that this product is not injective. 
For two adjacent essential vertices \(\sfv,\sfv'\) of \(\graf\), the edge \(\sfe\) between them splits into two edges in the decomposition. Hence their actions are not identical on the domain but identical on the codomain of the external product.

Let \(E^\sfv=\{\sfe^\sfv_1,\dots,\sfe^\sfv_{\ell(\sfv)+m(\sfv)}\}\) be the set of edges in \(\graf^\sfv\). For each triple \(1\le i<j<k\le 2\ell(\sfv)+m(\sfv)\) and each loop \(\sfe^\sfv_r\), which is the union of two half-edges \(\sfh^\sfv_i\) and \(\sfh^\sfv_{i+1}\), we define cycles in \(S_\bullet(\graf^\sfv)\)
\begin{align*}
\alpha^\sfv_{ijk}&\coloneqq \sfe(\sfh^\sfv_i)\cdot(\sfh^\sfv_j-\sfh^\sfv_k)+\sfe(\sfh^\sfv_j)\cdot(\sfh^\sfv_k-\sfh^\sfv_i)+\sfe(\sfh^\sfv_k)\cdot(\sfh^\sfv_i-\sfh^\sfv_j),\\
\beta^\sfv_r&\coloneqq \sfh^\sfv_{i+1}-\sfh^\sfv_{i},
\end{align*}
whose homology classes are called the \emph{star} and \emph{loop} classes in \(\graf^\sfv\).

\begin{definition}[Star-Loop submodule -- general case]\label{definition:star-loop basis}
We define \(\ring\)-submodules \(\Star(\graf^\sfv)\) and \(\Loop(\graf^\sfv)\) of \(H_1(B\graf^\sfv)\) as follows: for monomials \(p,q\in \ring[E^\sfv]\),
\begin{enumerate}
\item \(\Star(\graf^\sfv)\) is spanned over \(\ring\) by stabilized star classes \(p\cdot\alpha_{1jk}\) such that 
\begin{enumerate}
\item both \(\sfh^\sfv_j\) and \(\sfh^\sfv_k\) are the first half-edges of \(\sfe^\sfv_{j'}\) and \(\sfe^\sfv_{k'}\), respectively,
\item \(p\) is \(\sfe^\sfv_1\)-free, and \(\sfe^\sfv_{i'}\)-free for every \(j'<i'<k'\).
\end{enumerate}
\item \(\Loop(\graf^\sfv)\) is spanned over \(\ring\) by stabilized loop classes \(q\cdot\beta_r\) such that \(q\) is \(\sfe^\sfv_1\)-free.
\end{enumerate}
\end{definition}

Now we consider the sporadic case when the stem \(\sfT\) of \(\graf\cong(\sfT,\ell)\) has no edges. Namely, \(\graf\cong \graf_{\ell,0}\) is a bouquet of \(\ell\)-circles.
Then there is an embedding \(\iota:\graf_{\ell-1,2}\to\graf_{\ell,0}\), which induces a homomorphism \(\iota_*:H_\bullet(B\graf_{\ell-1,2})\to H_\bullet(B\graf_{\ell,0})\) between \(\ring[E]\)-modules.

Let \(\sfv\) be the central vertex of \(\graf_{\ell,0}\) and \(\{\sfe_1,\dots,\sfe_\ell\}\) denote edges of \(\graf_{\ell,0}\). Obviously, \(\graf_{\ell,0}^\sfv = \graf_{\ell,0}\).
\begin{definition}[Star-Loop submodule -- sporadic case]
We define a \(\ring[E]\)-submodule \(\SL(\graf_{\ell,0})\) as follows:
\begin{enumerate}
\item \(\Star(\graf_{\ell,0})\) is defined to be \(\iota_*(\Star(\graf_{\ell-2,2}))\).
\item \(\Loop(\graf_{\ell,0})\) is spanned over \(\ring\) by stabilized loop classes \(q\cdot\beta_r\) such that \(q\) is 
\begin{enumerate}
\item \(\sfe_{r-1}\)-free if \(1<r\le \ell\),
\item \(\sfe_\ell\)-free if \(1=r\).
\end{enumerate}
\end{enumerate}
\end{definition}

Observe that for a monomial \(p\), which is \(\sfe\)-free, \(p'\) is \(\sfe\)-free as well for any monomial \(p'\le p\).
Therefore we have the following obvious observation, whose proof will be omitted.
\begin{lemma}\label{lemma:subword action}
Suppose that \(p\cdot\alpha_{1jk}^\sfv\in \Star(\graf^\sfv)\) and \(q\cdot\beta_r^\sfv\in \Loop(\graf^\sfv)\) for some monomials \(p,q\in \ring[E^\sfv]\).
Then \(p'\cdot\alpha_{1jk}\in \Star(\graf^\sfv)\) and \(q'\cdot\beta_r\in \Loop(\graf^\sfv)\) for any monomials \(p'\le p\) and \(q'\le q\).
\end{lemma}

\begin{definition}[Star-Loop submodule]
We define a \(\ring\)-module \(\SL_\sfv\) and \(\SL(\graf)\) as 
\begin{align*}
\SL_\sfv&\coloneqq \ring\oplus\Star(\graf^\sfv)\oplus \Loop(\graf^\sfv)&
&\text{ and }&
\SL(\graf)&\coloneqq H_0(B\sfe_0)\otimes_\ring \left(\bigotimes_{\sfv\in V^\ess}\SL_\sfv\right).
\end{align*}
\end{definition}

\begin{theorem}\cite{AK2025}\label{theorem:basis}
The external product
\[
*: \SL(\graf)\to H_\bullet(B\graf)
\]
is an isomorphism between \(\ring\)-modules.
\end{theorem}

By using this isomorphism, we take a pull-back of the comultiplication \(\Sha^*_H\) on \(H_\bullet B\graf\) to \(\SL(\graf)\), which will be denoted by \(\Sha^*_H\) again.
\begin{remark}
Similar to \(S_\bullet(\graf)\), we may suppress "\(\otimes\)" to denote elements in \(\SL(\graf)\).
\end{remark}

Then \(H_0(B\sfe)\cong\ring[\sfe]\) as seen in \Cref{example:trivial}, and therefore for each \(a\ge 0\) and \(\gamma^\sfv\in \SL_\sfv\), we have
\begin{align*}
\Sha^*_H\left(\sfe_0^a\ast\left(\conv_{\sfv\in V^\ess}\gamma^\sfv\right)\right)&=
(-1)^s\sha(\sfe_0^a)\conv_{\sfv\in V^\ess}\Sha^*(\gamma^\sfv)
\end{align*}
for some \(s\in\bbZ\) by \Cref{proposition:commutative diagram with external product}.

\subsection{The formality}
In this subsection, we will consider the formality of \(S_\bullet(\graf)\). Roughly speaking, it is about the difference between \(S_\bullet(\graf)\) and its homology groups \(H_\bullet(B\graf)\) up to quasi-isomorphisms as \DGCA's.

We first focus on the elementary bunch of grapes \(\graf_{\ell,m}\) for \(m\ge1\).
\begin{proposition}\label{proposition:sha on basis}
For each \(\alpha_{ijk}\) and \(\beta_r\) in \(S_\bullet(\graf_{\ell,m})\),
\begin{align*}
\Sha^*_S(\alpha_{ijk}) &= \alpha_{ijk}\otimes 1+1\otimes\alpha_{ijk} +(\partial\otimes 1+1\otimes\partial) \mathsf{H}_{ijk},\\
\Sha^*_S(\beta_r)&= \beta_r\otimes 1+1\otimes \beta_r
\end{align*}
where 
\begin{align*}
\mathsf{H}_{ijk}&\coloneqq \sfh_i\otimes(\sfh_j-\sfh_k)+\sfh_j\otimes(\sfh_k-\sfh_i)+\sfh_k\otimes(\sfh_i-\sfh_j)\\
&=(\sfh_k-\sfh_j)\otimes \sfh_i+(\sfh_i-\sfh_k)\otimes \sfh_j+(\sfh_j-\sfh_i)\otimes \sfh_k \in S_1(\graf_{\ell,m})\otimes S_1(\graf_{\ell,m}).
\end{align*}

In particular, as elements in \(H(S_\bullet(\graf_{\ell,m})\otimes S_\bullet(\graf_{\ell,m}))\cong H(S_\bullet(\graf_{\ell,m}))\otimes H(S_\bullet(\graf_{\ell,m}))\),
\begin{align*}
[\Sha^*_S(\alpha_{ijk})]&=\alpha_{ijk}\otimes1+1\otimes\alpha_{ijk},& &\text{and}&
[\Sha^*_S(\beta_r)]&=\beta_r\otimes1+1\otimes\beta_r,
\end{align*}
which are primitive.
\end{proposition}
\begin{proof}
We have obviously that 
\[
\Sha^*_S(\beta_r)=\Sha^*_S(\sfh_{i+1}-\sfh_i)
=(\sfh_{i+1}\otimes1+1\otimes\sfh_{i+1})-(\sfh_i\otimes1+1\otimes\sfh_i)
=\beta_r\otimes1+1\otimes\beta_r.
\]

On the other hand,
\begin{align*}
\Sha^*_S(\alpha_{ijk})&= (\sfe(\sfh_i)\otimes 1+1\otimes \sfe(\sfh_i))\cdot\left( (\sfh_j-\sfh_k)\otimes 1 + 1\otimes(\sfh_j-\sfh_k)\right)\\
&\mathrel{\hphantom{=}}+(\sfe(\sfh_j)\otimes 1+1\otimes \sfe(\sfh_j))\cdot\left( (\sfh_k-\sfh_i)\otimes 1 + 1\otimes(\sfh_k-\sfh_i)\right)\\
&\mathrel{\hphantom{=}}+(\sfe(\sfh_k)\otimes 1+1\otimes \sfe(\sfh_k))\cdot\left( (\sfh_i-\sfh_j)\otimes 1 + 1\otimes(\sfh_i-\sfh_j)\right)\\
&=\alpha_{ijk}\otimes 1+1\otimes\alpha_{ijk}+(\partial\otimes 1+1\otimes\partial) \mathsf{H}_{ijk}.
\end{align*}
Here, we use the following fact 
\[
(\partial\otimes 1+1\otimes\partial)(\sfh\otimes\sfh')
= \partial\sfh\otimes\sfh' +(-1)^{\deg(\sfh)}\sfh\otimes\partial\sfh'
= \partial\sfh\otimes\sfh' -\sfh\otimes\partial\sfh'.
\]

Since \(H(S_\bullet(\graf_{\ell,m}))=H_0(S_\bullet(\graf_{\ell,m}))\oplus H_1(S_\bullet(\graf_{\ell,m}))\) is torsion-free over \(\ring\) by \Cref{theorem:basis for elementary}, we may identify \(\ring\)-modules
\[
H(S_\bullet(\graf_{\ell,m})\otimes S_\bullet(\graf_{\ell,m}))\cong H(S_\bullet(\graf_{\ell,m}))\otimes H(S_\bullet(\graf_{\ell,m})),
\]
which completes the proof.
\end{proof}

\begin{corollary}
Let \(p, q\in \ring[E]\). Then for \(\alpha_{ijk},\beta_r\in H(S_\bullet(\graf_{\ell,m}))\),
\begin{align*}
[\Sha^*_S(p\cdot\alpha_{ijk})]&=\sha^*(p)\cdot\left(\alpha_{ijk}\otimes1+1\otimes\alpha_{ijk}\right),\\
[\Sha^*_S(q\cdot\beta_{r})]&=\sha^*(q)\cdot\left(\beta_{r}\otimes1+1\otimes\beta_{r}\right)
\end{align*}
as elements in \(H(S_\bullet(\graf_{\ell,m}))\otimes H(S_\bullet(\graf_{\ell,m}))\).
\end{corollary}

Let \(A_\bullet(\graf_{\ell,m})\) be an acyclic chain complex over \(\ring[E]\) generated by \(a_{ijk}\) of degree \(0\) and \(\bar a_{ijk}\) of degree \(1\) for each \(1\le i<j<k\le 2\ell+m\), whose differential \(\partial_A\) is defined as
\begin{align*}
\partial_A(p(E)\cdot \bar a_{ijk}) = p(E)\cdot a_{ijk}
\end{align*}
for any \(p(E)\in \ring[E]\).
We consider a chain map \(\Phi_{\ell,m}:A_\bullet(\graf_{\ell,m})\to S_\bullet(\graf_{\ell,m})\) defined by
\[
\Phi_{\ell,m}(\bar a_{ijk})=\alpha_{ijk}\quad\text{ and }\quad
\Phi_{\ell,m}(a_{ijk})=0.
\]
Then since \(A_\bullet(\graf_{\ell,m})\) is acyclic, the mapping cone \((\tilde S_\bullet(\graf_{\ell,m}),\partial_{\tilde S})\) of \(\Phi_{\ell,m}\) is quasi-isomorphic to \(S_\bullet(\graf_{\ell,m})\).
More precisely, 
\[
\tilde S_\bullet(\graf_{\ell,m})\coloneqq A_\bullet(\graf_{\ell,m})[1]\oplus S_\bullet(\graf_{\ell,m}) = \ring[E]\left\langle \bar a_{ijk}, a_{ijk} : 1\le i<j<k\le 2\ell+m\right\rangle\oplus S_\bullet(\graf_{\ell,m})
\]
and
\begin{align*}
&\begin{cases}|\bar a_{ijk}|\coloneqq2,\\
|a_{ijk}|\coloneqq1,
\end{cases}
&
\begin{cases}
\partial_{\tilde S} \bar a_{ijk}\coloneqq a_{ijk}-\alpha_{ijk},\\
\partial_{\tilde S} a_{ijk}\coloneqq0.
\end{cases}
\end{align*}

\begin{proposition}
The chain complex \(\tilde S_\bullet(\graf_{\ell,m})\) admits a {\DGCA} structure.
\end{proposition}
\begin{proof}
Let us define a comultiplication \(\Sha^*_{\tilde S}\) for \(\bar a_{ijk}\) and \(a_{ijk}\) as
\begin{align*}
\Sha^*_{\tilde S}(a_{ijk})&\coloneqq a_{ijk}\otimes \varnothing+\varnothing\otimes a_{ijk},& &\text{and}&
\Sha^*_{\tilde S}(\bar a_{ijk})&\coloneqq \bar a_{ijk}\otimes \varnothing+\varnothing\otimes \bar a_{ijk}+\mathsf{H}_{ijk}.
\end{align*}
Then it suffices to prove the cocommutativity, coassociativity of \(\Sha^*\) and that \(\partial_{\tilde S}\) is a coderivation.

It is easy to check that \(\Sha^*_{\tilde S}\) on these generators are cocommutative via the map \(x\otimes y\mapsto (-1)^{|x|\cdot|y|} y\otimes x\).
For the generator \(\bar a_{ijk}\),
\begin{align*}
(\Sha^*_{\tilde S}\otimes 1)(\Sha^*_{\tilde S}(\bar a_{ijk}))
&=\Sha^*_{\tilde S}(\bar a_{ijk})\otimes1+\Sha^*_{\tilde S}(1)\otimes \bar a_{ijk}\\
&\mathrel{\hphantom{=}}
+\Sha^*_{\tilde S}(\sfh_i)\otimes(\sfh_j-\sfh_k)+\Sha^*_{\tilde S}(\sfh_j)\otimes(\sfh_k-\sfh_i)+\Sha^*_{\tilde S}(\sfh_k)\otimes(\sfh_i-\sfh_j)\\
&=\left(\bar a_{ijk}\otimes1\otimes1+1\otimes \bar a_{ijk}\otimes1+\mathsf{H}_{ijk}\otimes1\right)+1\otimes 1\otimes \bar a_{ijk}\\
&\mathrel{\hphantom{=}}+\left(1\otimes \mathsf{H}_{ijk}+ \sfh_i\otimes1\otimes(\sfh_j-\sfh_k)+\sfh_j\otimes1\otimes(\sfh_k-\sfh_i)+\sfh_k\otimes1\otimes(\sfh_i-\sfh_j)\right)\\
&=\bar a_{ijk}\otimes1\otimes1+\left(1\otimes \bar a_{ijk}\otimes1+1\otimes 1\otimes \bar a_{ijk}+1\otimes \mathsf{H}_{ijk}\right)\\
&\mathrel{\hphantom{=}}+\left(\mathsf{H}_{ijk}\otimes1+ \sfh_i\otimes1\otimes(\sfh_j-\sfh_k)+\sfh_j\otimes1\otimes(\sfh_k-\sfh_i)+\sfh_k\otimes1\otimes(\sfh_i-\sfh_j)\right)\\
&=\bar a_{ijk}\otimes\Sha^*(1)+1\otimes\Sha^*(\bar a_{ijk})\\
&\mathrel{\hphantom{=}}+(\sfh_k-\sfh_j)\otimes\Sha^*(\sfh_i)+(\sfh_i-\sfh_k)\otimes\Sha^*(\sfh_k)+(\sfh_j-\sfh_i)\otimes\Sha^*(\sfh_k)\\
&=(1\otimes\Sha^*_{\tilde S})(\Sha^*_{\tilde S}(\bar a_{ijk})),
\end{align*}
\begin{align*}
\Sha^*_{\tilde S}(\partial \bar a_{ijk}) &= \Sha^*_{\tilde S}(a_{ijk}-\alpha_{ijk})\\
&=\left(a_{ijk}\otimes\varnothing+\varnothing\otimes a_{ijk}\right)
-\left(\alpha_{ijk}\otimes\varnothing+\varnothing\otimes\alpha_{ijk}+(\partial\otimes 1+1\otimes\partial)\mathsf{H}_{ijk}\right)\\
&=(a_{ijk}-\alpha_{ijk})\otimes\varnothing+\varnothing\otimes(a_{ijk}-\alpha_{ijk})+(\partial_{\tilde S}\otimes 1+1\otimes\partial_{\tilde S})\mathsf{H}_{ijk}\\
&=\partial_{\tilde S} \bar a_{ijk}\otimes \varnothing+\varnothing\otimes\partial_{\tilde S} \bar a_{ijk}+(\partial_{\tilde S}\otimes 1+1\otimes\partial_{\tilde S})\mathsf{H}_{ijk}\\
&=(\partial_{\tilde S}\otimes 1+1\otimes\partial_{\tilde S})(\bar a_{ijk}\otimes \varnothing + \varnothing\otimes \bar a_{ijk} + \mathsf{H}_{ijk})\\
&=(\partial_{\tilde S}\otimes 1+1\otimes\partial_{\tilde S})\Sha^*_{\tilde S}(\bar a_{ijk}).\qedhere 
\end{align*}

For the generator \(a_{ijk}\), the coassociativity for \(a_{ijk}\) is obvious, and we have 
\[
(\partial_{\tilde S}\otimes 1+1\otimes\partial_{\tilde S})\Sha^*_{\tilde S}(a_{ijk}) = 0 = \Sha^*_{\tilde S}(\partial_{\tilde S} a_{ijk}).
\]

The counit \(e_{\tilde S}\) and coaugmentation \(u_{\tilde S}\) will be defined by pre- and post-compositions of \(e_S:S_\bullet(\graf_{\ell,m})\to \ring\) and \(u_S:\ring\to S_\bullet(\graf_{\ell,m})\) with the canonical quotient \(\tilde S_\bullet(\graf_{\ell,m})\to S_\bullet(\graf_{\ell,m})\) and inclusion \(S_\bullet(\graf_{\ell,m})\to \tilde S_\bullet(\graf_{\ell,m})\).
\end{proof}

\begin{proposition}\label{proposition:elementary formal}
The {\DGCA} \(S_\bullet(\graf_{\ell,m})\) is formal, and so is \(C_\bullet(B\graf_{\ell,m})\).
\end{proposition}
\begin{proof}
We have a sequence of quasi-isomorphisms between DG-coalgebras
\[
\begin{tikzcd}
C_\bullet(B\graf_{\ell,m})\ar[r,"g"] & S_\bullet(\graf_{\ell,m})\ar[r,hookrightarrow] & \tilde S_\bullet(\graf_{\ell,m}),
\end{tikzcd}
\]
where the right arrow is a canonical inclusion.
Hence it suffices to prove the existence of a quasi-isomorphism between DG-coalgebras \(H(B\graf_{\ell,m})\to \tilde S_\bullet(\graf_{\ell,m})\).

By \Cref{theorem:basis} and the followed discussion, there is an isomorphism between DG-coalgebras \(*:\SL(\graf_{\ell,m})\to H_\bullet(B\graf_{\ell,m})\).
Then we define a map \(\SL(\graf_{\ell,m})\to \tilde S_\bullet(\graf_{\ell,m})\) as \(p\alpha_{1jk}\mapsto p a_{1jk}\) and \(q\beta_r\mapsto q\beta_r\).
Obviously, this map induces an isomorphism on homology.

Moreover, this map preserves the coproduct by definitions of coproducts on \(\SL(\graf_{\ell,m})\) and \(\operatorname{Cone}(\Phi_{\ell,m})\). Therefore the zig-zag
\[
\begin{tikzcd}
C_\bullet(B\graf_{\ell,m})\ar[r,"g"] & S_\bullet(\graf_{\ell,m})\ar[r,hookrightarrow] & \tilde S_\bullet(\graf_{\ell,m})\ar[from=r]& \SL(\graf_{\ell,m}))\ar[r,"*"]& H(B\graf_{\ell,m})
\end{tikzcd}
\]
of quasi-isomorphisms of DG-coalgebras exist, and so both \(S_\bullet(\graf_{\ell,m})\) and \(C_\bullet(B\graf_{\ell,m})\) are formal as {\DGCA}'s over \(\ring\).
\end{proof}

Let \(\graf=(V,E)\cong(\sfT,\ell)\) be a bunch of grapes.
Then as seen earlier, there is a unique decomposition of \(\graf\) into elementary bunches of grapes \(\{\graf^\sfv:\sfv\in V^{\ess}\}\).

For each \(\sfv\in V^\ess\), we consider an acyclic chain complex \(A_\bullet(\graf^\sfv)\cong A_\bullet(\graf_{\ell(\sfv),m(\sfv)})\) and the chain map \(\Phi^\sfv:A_\bullet(\graf^\sfv)\to S_\bullet(\graf^\sfv)\) between \(\ring[E^\sfv]\)-modules.

Now we extend the ring by tensoring \(\ring[E]\).
That is, 
\[
\Phi^\sfv_E:\ring[E]\otimes_{\ring[E^\sfv]} A_\bullet(\graf^\sfv) \to \ring[E]\otimes_{\ring[E^\sfv]} S_\bullet(\graf^\sfv).
\]
We denote the mapping cone of \(\Phi^\sfv_E\) by \(\tilde S_\bullet^E(\graf^\sfv)\).

Then we define a new chain complex 
\[
\tilde S_\bullet(\graf)\coloneqq\bigwedge_{\sfv\in V^\ess} \tilde S_\bullet^E(\graf^\sfv)
=\bigwedge_{\sfv\in V^\ess} \ring[E]\otimes_{\ring[E^\sfv]}\tilde S_\bullet(\graf^\sfv)
=\bigwedge_{\sfv\in V^\ess} \ring[E]\otimes_{\ring[E^\sfv]}\left(A_\bullet^\sfv[1]\oplus S_\bullet(\graf^\sfv)\right),
\]
which admits a {\DGCA} structure with respect to \(\Sha^*_{\tilde S}\) defined as follows:
For \(p\in \ring[E]\) and \(\gamma^\sfv\in \tilde S_\bullet(\graf^\sfv)\),
\[
\Sha^*_{\tilde S}\left(p\cdot \bigwedge_{\sfv\in V^\ess} \gamma^\sfv\right) = (-1)^s \sha^*(p)\cdot\left( \conv_{\sfv\in V^\ess} \Sha^*_{\tilde S}(\gamma^\sfv)\right)
\]
for some \(s\in\bbZ\) by \Cref{proposition:commutative diagram with external product} as before.

\begin{proposition}
There is an inclusion \(S_\bullet(\graf)\to \tilde S_\bullet(\graf)\), which is a quasi-isomorphism between {\DGCA}'s.
\end{proposition}
\begin{proof}
It is obvious from definition that the canonical inclusion 
\[
\bigotimes_{\sfv\in V^\ess} \ring[E]\otimes_{\ring[E^\sfv]}S_\bullet(\graf^\sfv) \to \tilde S_\bullet(\graf)
\]
is a quasi-isomorphism between {\DGCA}'s.
Moreover, since \(\ring[E]\otimes_{\ring[E^\sfv]}S_\bullet(\graf^\sfv)\cong S_\sfv\), and so we have an isomorphism 
\[
S_\bullet(\graf)=\bigwedge_{\sfv\in V^\ess} S_\sfv\cong \bigwedge_{\sfv\in V^\ess} \ring[E]\otimes_{\ring[E^\sfv]}S_\bullet(\graf^\sfv)
\]
between {\DGCA}'s as well. 
\end{proof}

\begin{theorem}\label{theorem:formality}
Let \(\graf=(V,E)\cong(\sfT,\ell)\) be a bunch of grapes.
Then the {\DGCA} \(S_\bullet(\graf)\) is formal, and so is \(C_\bullet(\graf)\).
\end{theorem}
\begin{proof}
We will use essentially the same argument as the proof of \Cref{proposition:elementary formal}.

As before, we already have quasi-isomorphisms between DG-coalgebras
\[
\begin{tikzcd}
C_\bullet(B\graf)\ar[r,"g"] & S_\bullet(\graf)\ar[r,hookrightarrow] & \tilde S_\bullet(\graf)
\end{tikzcd}\quad\text{ and }\quad
\begin{tikzcd}
\SL(\graf)\ar[r,"*"]&H_\bullet(B\graf).
\end{tikzcd}
\]
Therefore it suffices to prove the existence of a quasi-isomorphism between DG-coalgebras \(\SL(\graf)\to \tilde S_\bullet(\graf)\).

We define a map \(\SL(\graf)\to \tilde S_\bullet(\graf)\) as \(p\alpha_{1jk}\mapsto p a_{1jk}\) and \(q\beta_r\mapsto q\beta_r\).
Then it is easy to check that this map is compatible with comultiplications \(\Sha^*_{\SL}\) and \(\Sha^*_{\tilde S}\) and induces an isomorphism between {\DGCA}'s, which completes the proof.
\end{proof}

\subsection{The primitivity}
Recall \Cref{definition:R0} for the subring \(\ring_0[E]\) of \(\ring[E]\) generated by elements of the form \(\sfe_i-\sfe_j\) over \(\ring\).
Then we prove the following theorem about the primitivity of a homology class in \(H_\bullet(B\graf)\).
\begin{theorem}\label{theorem:primitivity}
Let \(\graf=(V,E)\cong(\sfT,\ell)\) be a bunch of grapes.
Then \(\alpha\in H_\bullet(B\graf)\) is primitive with respect to the comultiplication \(\Sha^*_H\) if and only if \(\alpha\) is a linear combination over \(\ring\) of the following: for a fixed edge \(\sfe\in E\),
\begin{enumerate}
\item \(\sfe\cdot\varnothing\in H_0(B_0\graf)\),
\item \(p\alpha_{1jk}^\sfv\in \SL(\graf)\), or
\item \(q\beta_r^\sfv\in \SL(\graf)\),
\end{enumerate}
where \(p,q\in \ring_0[E]\) and \(\sfv\in V^\ess\).
\end{theorem}

We first consider the `if' part.
\begin{lemma}
Let \(p=(E^\sfv)^\bfa\in R[E^\sfv]\) be a monomial and \(\gamma^\sfv\in \SL_\sfv\) be either \(\alpha_{1jk}^\sfv\) or \(\beta_r^\sfv\). Then
\begin{align*}
\Sha^*_H((E^\sfv)^\bfa\gamma^\sfv) &=\sum_{\bfb\le\bfa}\binom{\bfa}{\bfb}\left((E^\sfv)^\bfb\gamma^\sfv\otimes \sfe_0^{\|\bfa-\bfb\|}+
\sfe_0^{\|\bfb\|}\otimes (E^\sfv)^{\bfa-\bfb}\cdot\gamma^\sfv\right),
\end{align*}
which is an element in \(\SL_\sfv\otimes \ring[\sfe_0]\oplus \ring[\sfe_0]\otimes\SL_\sfv\).
\end{lemma}
\begin{proof}
Since \(\Sha^*_H(\gamma^\sfv)=\gamma^\sfv\otimes1+1\otimes\gamma^\sfv\), and by \Cref{lemma:comultiplication on ring},
\begin{align*}
\Sha^*_H((E^\sfv)^\bfa\gamma^\sfv) &= \sha^*((E^\sfv)^\bfa)\Sha^*(\gamma^\sfv)\\
&=\sum_{\bfb\le\bfa}\binom{\bfa}{\bfb}((E^\sfv)^\bfb\otimes (E^\sfv)^{\bfa-\bfb})(\gamma^\sfv\otimes 1+1\otimes\gamma^\sfv)\\
&=\sum_{\bfb\le\bfa}\binom{\bfa}{\bfb}\left((E^\sfv)^\bfb\gamma^\sfv\otimes (E^\sfv)^{\bfa-\bfb}\cdot1+
(E^\sfv)^\bfb\cdot1\otimes (E^\sfv)^{\bfa-\bfb}\cdot\gamma^\sfv\right)\\
&=\sum_{\bfb\le\bfa}\binom{\bfa}{\bfb}\left((E^\sfv)^\bfb\gamma^\sfv\otimes \sfe_0^{\|\bfa-\bfb\|}+
\sfe_0^{\|\bfb\|}\otimes (E^\sfv)^{\bfa-\bfb}\cdot\gamma^\sfv\right).
\end{align*}

Finally, by \Cref{lemma:subword action}, \((E^\sfv)^\bfb\gamma^\sfv\in\SL_\sfv\) for every \(\bfb\le\bfa\) and we are done.
\end{proof}

\begin{lemma}\label{lemma:primitivity if}
Let \(p\in \ring_0[E]\) and \(\gamma^\sfv\in\SL_\sfv(\graf)\) be either \(\alpha_{1jk}^\sfv\) or \(\beta_r^\sfv\) for some vertex \(\sfv\in V^\ess\).
Then \(p\cdot\gamma^\sfv\) is primitive.
\end{lemma}
\begin{proof}
Without loss of generality, we may assume that \(p\) and \(q\) are homogeneous elements.
By definition, 
\begin{align*}
\Sha^*_H(p\cdot \gamma^\sfv)&=\sha^*(p)\cdot\Sha^*(\gamma^\sfv)
=\sha^*(p)\cdot(\gamma^\sfv\otimes\varnothing+\varnothing\otimes\gamma^\sfv).
\end{align*}

Let \(F:\ring[E]\otimes \ring[E]\to \ring[\sfe_0]\) for a fixed \(\sfe_0\in E\).
Then \(\sha^*(p)\cdot(\gamma^\sfv\otimes \varnothing)=F(\sha^*(p))\cdot(\gamma^\sfv\otimes \varnothing)\) since \(\sfe\cdot\varnothing=\sfe_0\cdot\varnothing\) for any \(\sfe\in E\).
Hence it is the same as \((p\cdot\gamma^\sfv)\otimes \varnothing\) by \Cref{proposition:graded augmentation}, and similarly, \(\sha^*(p)\cdot(1\otimes\gamma^\sfv)=1\otimes(p\cdot\gamma^\sfv)\).
Therefore \(p\cdot\gamma^\sfv\) is primitive.
\end{proof}

\begin{proof}[Proof of \Cref{theorem:primitivity}]
Obviously, 
\[
\Sha^*_H(\sfe\cdot\varnothing) = \sha^*(\sfe)\cdot\Sha^*_H(\varnothing)
=(\sfe\otimes 1+1\otimes \sfe)\cdot(\varnothing\otimes\varnothing)
=(\sfe\cdot\varnothing)\otimes 1 + 1\otimes (\sfe\cdot\varnothing).
\]
Hence \(\sfe\cdot\varnothing\) is primitive.
Therefore since \(\Sha_H^*\) is \(\ring\)-linear together with \Cref{lemma:primitivity if}, the `if' part is done.

Suppose that \(c\) is primitive. That is, \(\Sha^*(c)=c\otimes 1+1\otimes c\).
Since \(\Sha^*\) preserves the total bigrading, we may assume that \(c\in H_i(B_n\graf_{\ell,m})\) for some \(i,n\ge0\).

If \(i=n=0\), then \(c\) is a counit and so \(c=u\cdot\varnothing\) for some \(u\in\ring\). 
Therefore \(\Sha^*(c)=u\cdot\varnothing\otimes\varnothing\neq c\otimes 1+1\otimes c\). This is a contradiction and \(i+n\ge1\).

Suppose that \(i=0\), or equivalently, \(c\in H_0(B_n\graf)\cong\ring[E]\) for some \(n\ge 1\). Then since \(H_0(B_n\graf_{\ell,m})\cong\ring\langle \sfe^n\cdot\varnothing\rangle\) for some \(\sfe\in E\), we have \(c=[u\sfe^n\cdot\varnothing]\) for some \(0\neq u\in\ring\). However,
\begin{align*}
\Sha^*(c)&=\Sha^*([u\sfe^n\cdot\varnothing])\\
&=u[(\sha^*(\sfe)\cdots\sha^*(\sfe))\cdot\Sha^*(\varnothing)]\\
&=u[(\sfe\otimes1+1\otimes\sfe)\cdots(\sfe\otimes1+1\otimes\sfe)\cdot(\varnothing\otimes\varnothing)]\\
&=u\sum_{i=0}^n \binom ni \left[\sfe^i\cdot\varnothing\otimes \sfe^{n-i}\cdot\varnothing\right],
\end{align*}
which is primitive if and only if \(n=1\).

For \(\i\ge 1\), by \Cref{theorem:basis}, we may assume that \(c\) is an element in \(\SL\).
That is, \(c\) is the image under \(*\) of a linear combination over \(\ring\) of elements of the form
\[
\sfe_0^a\bigotimes_{\sfv\in V^\ess}p^\sfv\gamma^\sfv,
\]
where \(\gamma^\sfv\in \SL_\sfv\) is either \(1,\alpha^\sfv_{1jk}\), or \(\beta^\sfv_r\) for some \(p^\sfv\in \ring[E^\sfv]\), and integers \(j,k,r\).

Since each \(p^\sfv\) can be pulled out in \(H_\bullet(B\graf)\), \(c\) can be expressed as 
\[
c=\sum_s p_s c_s,\quad \text{for}\quad c_s=\conv_{\sfv\in V^\ess} \gamma_s^\sfv.
\]
for some polynomial \(p_s\in \ring[E]\). 
We may assume that sequences \((\gamma_s^\sfv)_{\sfv\in V^\ess}\) are pairwise different with respect to \(s\) by collecting terms further if necessary, which implies that the external products \(c_s\)'s are linearly independent by \Cref{lemma:subword action}.

Then the comultiplication \(\Sha_H^*(c)\) is as follows: by \Cref{proposition:commutative diagram with external product},
\begin{align*}
\Sha^*_H(c)&=
\sum_s\Sha^*_H\left(
p_s\conv_{\sfv\in V^\ess} \gamma_s^\sfv
\right)\\
&=\sum_s \sha^*(p_s) \left(\conv_{\sfv\in V^\ess}\otimes\conv_{\sfv\in V^\ess}\right) T\left(
\bigotimes_{\sfv\in V^\ess} \Sha^*_H(\gamma_s^\sfv)
\right)\\
&=\sum_s \sha^*(p_s) \left(\conv_{\sfv\in V^\ess}\otimes\conv_{\sfv\in V^\ess}\right) T\left(
\bigotimes_{\sfv\in V^\ess} \left(\gamma_s^\sfv\otimes1+1\otimes \gamma_s^\sfv\right)
\right).
\end{align*}
Hence each summand of \(\Sha^*_H(c)\) up to sign is of the form 
\[
q_s \left(\conv_{\sfv\in V_1^\ess}\gamma_s^\sfv\right)\otimes\left(\conv_{\sfv\in V_2^\ess}\gamma_s^\sfv\right)
\]
for some subsets \(V_1^\ess, V_2^\ess\) of \(V^\ess\).
This is nothing but a separation of \(\gamma^\sfv_s\)'s, and so summands in \(\Sha^*_H(c)\) with respect to \(s\) are linearly independent since so are \(c_s\)'s as observed above.
Hence \(c\) is primitive if and only if so is \(p_s c_s\) for each \(s\), and we may assume that 
\[
c= p \conv_{\sfv\in V^\ess} \gamma^\sfv = \sfe_0^a \conv_{\sfv\in V^\ess} p^\sfv\gamma^\sfv.
\]
by pushing each \(p^\sfv\) back inside the external product.

Now let us take a look at \(\Sha^*_H(c)\) again. Then for a chosen vertex \(\sfw\in V^\ess\) with \(\gamma^\sfw\neq 0\), namely, \(\gamma^\sfw\) is either \(\alpha^\sfw_{1jk}\) or \(\beta^\sfw_r\), we collect all terms whose first factor is \(\gamma^\sfw\).
Then it is given as
\[
\gamma^\sfw \otimes \left(
\sfe_0^a \conv_{\substack{\sfv\in V^\ess\\\sfv\neq\sfw}} p^\sfw\gamma^\sfv
\right).
\]
Hence \(c\) is primitive only if the second factor of this summand vanishes, which is equivalent to all \(\gamma^\sfv\) except for \(\gamma^\sfw\) must be trivial and so \(c\) is of the form \(\sfe_0^a p^\sfw \gamma^\sfw\).

Finally, for \(p=\sfe_0^a p^\sfw\), we have 
\begin{align*}
\Sha^*_H(c) &= \sha^*(p)(\gamma^\sfw\otimes 1+1\otimes\gamma^\sfw)
=\sha^*(p)(\gamma^\sfw\otimes 1)+\sha^*(p)(1\otimes\gamma^\sfw)
=p\gamma^\sfw\otimes1 + 1\otimes p\gamma^\sfw
\end{align*}
since \(c\) is primitive.
Then \(p\in \ring_0[E]\) by \Cref{proposition:graded augmentation} and we are done.
\end{proof}

\bibliographystyle{alpha}
\bibliography{references.bib}

@misc{AK2025,
      title={Hilbert polynomials of configuration spaces over graphs of circumference at most 1}, 
      author={Byung Hee An and Jang Soo Kim},
      year={2025},
      eprint={2505.24416},
      archivePrefix={arXiv},
      primaryClass={math.GT},
      url={https://arxiv.org/abs/2505.24416}, 
}

@article {Ramos2018,
    AUTHOR = {Ramos, Eric},
     TITLE = {Stability phenomena in the homology of tree braid groups},
   JOURNAL = {Algebr. Geom. Topol.},
  FJOURNAL = {Algebraic \& Geometric Topology},
    VOLUME = {18},
      YEAR = {2018},
    NUMBER = {4},
     PAGES = {2305--2337},
      ISSN = {1472-2747,1472-2739},
   MRCLASS = {05C25 (05C05 05C10 05E40 57M15)},
  MRNUMBER = {3797068},
       DOI = {10.2140/agt.2018.18.2305},
       URL = {https://doi.org/10.2140/agt.2018.18.2305},
}

@article {ADCK2020,
    AUTHOR = {An, Byung Hee and Drummond-Cole, Gabriel C. and Knudsen, Ben},
     TITLE = {Edge stabilization in the homology of graph braid groups},
   JOURNAL = {Geom. Topol.},
  FJOURNAL = {Geometry \& Topology},
    VOLUME = {24},
      YEAR = {2020},
    NUMBER = {1},
     PAGES = {421--469},
      ISSN = {1465-3060,1364-0380},
   MRCLASS = {20F36 (05C40 55R80)},
  MRNUMBER = {4080487},
MRREVIEWER = {Mehmet\ Emin\ Akta\c s},
       DOI = {10.2140/gt.2020.24.421},
       URL = {https://doi.org/10.2140/gt.2020.24.421},
}

@article {KP2012,
    AUTHOR = {Ko, Ki Hyoung and Park, Hyo Won},
     TITLE = {Characteristics of graph braid groups},
   JOURNAL = {Discrete Comput. Geom.},
  FJOURNAL = {Discrete \& Computational Geometry. An International Journal
              of Mathematics and Computer Science},
    VOLUME = {48},
      YEAR = {2012},
    NUMBER = {4},
     PAGES = {915--963},
      ISSN = {0179-5376,1432-0444},
   MRCLASS = {20F36 (05C25)},
  MRNUMBER = {3000570},
MRREVIEWER = {J.\ S.\ Birman},
       DOI = {10.1007/s00454-012-9459-8},
       URL = {https://doi.org/10.1007/s00454-012-9459-8},
}

@incollection {Ghrist2001,
    AUTHOR = {Ghrist, Robert},
     TITLE = {Configuration spaces and braid groups on graphs in robotics},
 BOOKTITLE = {Knots, braids, and mapping class groups---papers dedicated to
              {J}oan {S}. {B}irman ({N}ew {Y}ork, 1998)},
    SERIES = {AMS/IP Stud. Adv. Math.},
    VOLUME = {24},
     PAGES = {29--40},
 PUBLISHER = {Amer. Math. Soc., Providence, RI},
      YEAR = {2001},
      ISBN = {0-8218-2966-1},
   MRCLASS = {55R80 (20F36 93C85)},
  MRNUMBER = {1873106},
MRREVIEWER = {Vagn\ Lundsgaard\ Hansen},
       DOI = {10.1090/amsip/024/03},
       URL = {https://doi.org/10.1090/amsip/024/03},
}

@book {Abrams2000,
    AUTHOR = {Abrams, Aaron David},
     TITLE = {Configuration spaces and braid groups of graphs},
      NOTE = {Thesis (Ph.D.)--University of California, Berkeley},
 PUBLISHER = {ProQuest LLC, Ann Arbor, MI},
      YEAR = {2000},
     PAGES = {67},
      ISBN = {978-0599-85818-3},
   MRCLASS = {99-05},
  MRNUMBER = {2701024},
       URL =
              {http://gateway.proquest.com/openurl?url_ver=Z39.88-2004&rft_val_fmt=info:ofi/fmt:kev:mtx:dissertation&res_dat=xri:pqdiss&rft_dat=xri:pqdiss:9979537},
}
\end{document}